\documentclass{article}
\usepackage[french,british]{babel}

\usepackage[T1]{fontenc}

\usepackage{amsmath}
\usepackage{amsfonts}
\usepackage{amsbsy}
\usepackage{amssymb}

\usepackage{amsthm}

\theoremstyle{plain}

\newtheorem{thm}{Theorem}[section]
\newtheorem{lem}{Lemma}[section]
\newtheorem{nas}{Corrolary}[section]

\newtheorem{ozn}{Definition}[section]
\newtheorem{exm}{Example}

\begin{document}

\selectlanguage{british}

\title{On minimax estimation problem for stationary stochastic sequences from observations in special sets of points}

\author{
Oleksandr Masyutka\thanks{Department of Mathematics and Theoretical Radiophysics, Taras Shevchenko National University of Kyiv, 01601, Ukraine, e-mail: omasyutka@gmail.com},
Mikhail Moklyachuk\thanks
{Department of Probability Theory, Statistics and Actuarial
Mathematics, Taras Shevchenko National University of Kyiv, Kyiv 01601, Ukraine, e-mail: moklyachuk@gmail.com}
}

\date{\today}

\maketitle

\renewcommand{\abstractname}{Abstract}
\begin{abstract}
The problem  of the mean-square optimal estimation of the linear functionals  which depend on the unknown values of a stochastic stationary sequence from observations of the sequence in special sets of points is considered. Formulas for calculating the mean-square error and the spectral characteristic of the optimal linear estimate of the functionals are derived under the condition of spectral certainty, where the spectral density of the sequence is exactly known. The minimax (robust)
method of estimation is applied in the case where the spectral density of the sequence is not known exactly while some sets of admissible spectral densities are given. Formulas that determine the least favourable spectral densities and the minimax spectral characteristics are derived for some special sets of admissible densities.
\end{abstract}

\vspace{2ex}

\maketitle

\vspace{2ex}

\section{Introduction}

The problem of estimation of the unknown values of stochastic processes is of constant interest in the theory of stochastic processes. The formulation of the interpolation, extrapolation and filtering problems for stationary stochastic sequences with known spectral densities and reducing them to the corresponding problems of the theory of functions belongs to Kolmogorov [13]. Effective methods of solution of the estimation problems for stationary stochastic sequences and processes were developed by Wiener [29] and Yaglom [30]. Further results are presented in the books  by Rozanov [24] and Hannan [8].
The crucial assumption of most of the methods developed for estimating the unobserved values of stochastic processes is that the spectral densities of the involved stochastic processes are exactly known. However, in practice complete information on the spectral densities is impossible in most cases.
In this situation one finds parametric or nonparametric estimate of the unknown spectral density and then apply one of the traditional estimation methods provided that the selected density is the true one. This procedure can result in significant increasing of the value of error as Vastola and Poor [28] have demonstrated with the help of some examples.
To avoid this effect one can search the estimates which are optimal for all densities from a certain class of admissible spectral densities. These estimates are called minimax since they minimize the maximum value of the error.
The paper by Grenander [7] was the first one where this approach to extrapolation problem for stationary processes was proposed.
Several models of spectral uncertainty and minimax-robust methods of data processing can be found in the survey paper by  Kassam and  Poor [12]. Franke [3], Franke and Poor [4] investigated the minimax extrapolation  and filtering problems for stationary sequences with the help of convex optimization methods. This approach makes it possible to find equations that determine the least favorable spectral densities for different classes of densities.
In the papers by Moklyachuk {17. 18] and  Liu, Xue and Taniguchi [15]
the problems of extrapolation, interpolation and filtering  for  functionals which depend on the unknown values of stationary processes and sequences are investigated.
The estimation problems for functionals which depend on the unknown values of multivariate stationary stochastic processes is the aim of the book by  Moklyachuk and Masytka [19],
Golichenko and Moklyachuk [6] investigated the interpolation, extrapolation and filtering problems for periodically correlated stochastic sequences.
In their book Luz and Moklyachuk [16]  deals with the estimation problems for functionals which depend on the unknown values of stochastic sequences with stationary increments.
Prediction of stationary processes with missing observations  is investigated in papers by Bondon [1, 2], Kasahara, Pourahmadi and Inoue [11, 22].
Moklyachuk and Sidei [26,27]
investigated problems of the interpolation, extrapolation and filtering of stationary stochastic sequences and processes with missing observations. See also the book by Moklyachuk, Masytka and Sidei [20].

In this article we consider the problem of the mean-square optimal estimation of the linear
functionals
\[A_{S_1}\xi=\sum\limits_{j=-\infty}^{-M_1-1}a(j)\xi(j)+\sum\limits_{j=0}^{N}a(j)\xi(j),\; A_{S_2}\xi=\sum\limits_{j=0}^{N}a(j)\xi(j)+\sum\limits_{j=N+M_2+1}^{\infty}a(j)\xi(j),\]
\[A_{S_3}\xi=\sum\limits_{j=-\infty}^{-M_1-1}a(j)\xi(j)+
\sum\limits_{j=0}^{N}a(j)\xi(j)+\sum\limits_{j=N+M_2+1}^{\infty}a(j)\xi(j),\]
which depend on the unknown values of a stochastic stationary sequence $\xi(j),\,j\in \mathbb{Z}$  from observations of the sequence at points $j\in S_k$ respectively, where
\[S_1=\{-M_1, \ldots, -1\}\cup \{N+1, N+2, \ldots\}, \; S_2=\{\ldots, -2, -1\}\cup  \{N+1, \ldots, N+M_2\},\]
\[S_3=\{-M_1, \ldots, -1\}\cup \{N+1, \ldots, N+M_2\}.\]
The problem is investigated in the case of spectral certainty, where the spectral density of the sequence is exactly known, and in the case of spectral uncertainty, where the spectral density of the sequence is unknown while a class of admissible spectral densities is given.

\section{The classical Hilbert space projection method of linear estimation}

Let $\xi(j),\,j\in \mathbb{Z}$,  be a (wide sense) stationary stochastic sequence. We will consider $\xi(j)$ as elements of the Hilbert space  $H=L_2(\Omega,\mathcal{F},P)$ of complex valued random variables with zero first moment, $E\xi=0$, finite second moment, $E|\xi|^2<\infty$, and the inner product $\left(\xi,{\eta}\right)=E\xi\overline{\eta}$. The correlation function $r(n)=E\xi(j+n)\overline{\xi(j)}$ of the stationary stochastic sequence $\xi(j),\,j\in \mathbb{Z}$, admits the spectral decomposition [5]
\[
 r(n)=\int\limits_{-\pi}^{\pi}e^{in\lambda}F(d\lambda),
 \]
where $F(d\lambda)$ is the spectral measure of the sequence.
We will consider stationary stochastic sequences with absolutely continuous spectral measures and the correlation functions of the form
\[
 r(n)=\,\frac{1}{2\pi}\int\limits_{-\pi}^{\pi}e^{ik\lambda}f(\lambda)d\lambda,
 \]
where $f(\lambda)$ is the spectral density function of the sequence $\xi(j)$ that satisfies the minimality condition
\begin {equation} \label{min}
 \int_{- \pi}^{ \pi} \frac {1} {f( \lambda)} \,d \lambda< \infty .
\end{equation}

This condition is necessary and sufficient in order that the error-free interpolation of unknown values of the sequence is impossible [24].

The stationary stochastic sequence $\xi(j),\,j\in \mathbb{Z}$, admits the spectral decomposition [5, 10]
\begin {equation} \label{sp-r}
\xi(j)=\int\limits_{-\pi}^{\pi}e^{ij\lambda}{Z}(d\lambda),
\end{equation}
where ${Z}(\Delta)$ is the orthogonal stochastic measure of the sequence such that
\[
EZ(\Delta_1)\overline{Z(\Delta_2)}=F(\Delta_1\cap\Delta_2)=\,\frac{1}{2\pi}\int_{\Delta_1\cap\Delta_2}f(\lambda)d\lambda.
\]

Consider the problem of the mean-square optimal estimation of the linear functionals
\[A_{S_1}\xi=\sum\limits_{j=-\infty}^{-M_1-1}a(j)\xi(j)+\sum\limits_{j=0}^{N}a(j)\xi(j),\; A_{S_2}\xi=\sum\limits_{j=0}^{N}a(j)\xi(j)+\sum\limits_{j=N+M_2+1}^{\infty}a(j)\xi(j),\]
\[A_{S_3}\xi=\sum\limits_{j=-\infty}^{-M_1-1}a(j)\xi(j)+
\sum\limits_{j=0}^{N}a(j)\xi(j)+\sum\limits_{j=N+M_2+1}^{\infty}a(j)\xi(j),\]
which depend on the unknown values of a stochastic stationary sequence $\xi(j),j\in\mathbb{Z},$  from observations of the sequence at points $j\in S_k$ respectively, where
\[S_1=\{-M_1, \ldots, -1\}\cup \{N+1, N+2, \ldots\}, \; S_2=\{\ldots, -2, -1\}\cup  \{N+1, \ldots, N+M_2\},\]
\[S_3=\{-M_1, \ldots, -1\}\cup \{N+1, \ldots, N+M_2\}.\]

It follows from the spectral decomposition (2) of the sequence \(\xi(j)\) that we can represent the functionals $A_{S_1}\xi$ and $A_{S_2}\xi$ in the following form
\begin{equation}\label{asx1}
A_{S_1}\xi=\int\limits_{-\pi}^{\pi}\left(A_{-M_1}(e^{i\lambda})+A_N(e^{i\lambda})\right){Z}(d\lambda),\;
A_{S_2}\xi=\int\limits_{-\pi}^{\pi}\left(A_{N}(e^{i\lambda})+A_{N+M_2+}(e^{i\lambda})\right){Z}(d\lambda)
 \end{equation}
\begin{equation}\label{asx2}
A_{S_3}\xi=\int\limits_{-\pi}^{\pi}\left(A_{-M_1}(e^{i\lambda})+A_N(e^{i\lambda})+A_{N+M_2+}(e^{i\lambda})\right){Z}(d\lambda)
 \end{equation}
where
\begin{equation*}
A_{-M_1}(e^{i\lambda})=\sum\limits_{j=-\infty}^{-M_1-1}a(j)e^{ij\lambda}, \; A_N(e^{i\lambda})=\sum\limits_{j=0}^{N}a(j)e^{ij\lambda },\; A_{N+M_2+}(e^{i\lambda})=\sum\limits_{j=N+M_2+1}^{\infty}a(j)e^{ij\lambda }.
 \end{equation*}

Denote by $H^{S_k}(\xi)$ the subspaces  of the Hilbert space  $H=L_2(\Omega,\mathcal{F},P)$ generated by elements $\{\xi(j):j\in S_k\}$ respectively. Denote by $L_{2} (f)$ the Hilbert space of functions $a( \lambda )$ such that
 \[ \int_{- \pi}^{ \pi} \left|a(\lambda)\right|^2 f(\lambda) d \lambda  < \infty. \] Denote by $L_2^{S_k}(f)$ the subspaces of $L_2(f)$ generated by the functions $\{e^{in \lambda}, n \in S_k\}$ respectively.
The mean square optimal linear estimates $\hat{A}_{S_k}\xi$ of the functionals $A_{S_k}\xi$ from observations of the sequence \(\xi(j)\) at points of time  $j\in S_k$ are elements of the $H^{S_k}(\xi)$ respectively. They can be represented in the form
\begin{equation}\label{hatas}
\hat{A}_{S_k}\xi=\,\int\limits_{-\pi}^{\pi}h_k(e^{i\lambda}){Z}(d\lambda),
 \end{equation}
where $h_k(e^{i\lambda}) \in L_2^{S_k}(f)$ are the spectral characteristics of the estimates $\hat{A}_{S_k}\xi$ respectively.

The mean square errors $\Delta(h_k;f)=E\left|A_{S_k}\xi-\hat{A}_{S_k}\xi\right|^2$ of the estimates $\hat{A}_{S_k}\xi$ are given by the formulas
\begin{equation*}
\Delta(h_1;f)=\,\frac{1}{2\pi}\int\limits_{-\pi}^{\pi}\left|A_{-M_1}(e^{i\lambda})+A_N(e^{i \lambda})-h_1(e^{i \lambda})\right|^{2} f( \lambda ) d\lambda,
\end{equation*}
\begin{equation*}
\Delta(h_2;f)=\,\frac{1}{2\pi}\int\limits_{-\pi}^{\pi}\left|A_N(e^{i \lambda})+A_{N+M_2+}(e^{i \lambda})-h_2(e^{i \lambda})\right|^{2} f( \lambda ) d\lambda.
\end{equation*}
\begin{equation*}
\Delta(h_3;f)=\,\frac{1}{2\pi}\int\limits_{-\pi}^{\pi}\left|A_{-M_1}(e^{i\lambda})+A_N(e^{i \lambda})+A_{N+M_2+}(e^{i \lambda})-h_3(e^{i \lambda})\right|^{2} f( \lambda ) d\lambda,
\end{equation*}
The Hilbert space projection method proposed by  Kolmogorov [13] makes it possible to find the spectral characteristics  $h_k(e^{i\lambda})$ and the mean square errors $\Delta(h_k;f)$ of the optimal linear estimates of the functionals $A_{S_k}\xi$ respectively in the case where the spectral density $f(\lambda)$ of the sequence is exactly known and the minimality condition (1) is satisfied. The spectral characteristic can be found from the following conditions:
\begin{equation*} \begin{split}
1)&\; h_k(e^{i\lambda}) \in L_2^{S_k}(f);\\
2)&\; A_{-M_1}(e^{i\lambda})+A_N(e^{i\lambda})- h_1(e^{i\lambda}) \bot \, L_2^{S_1}(f),\; A_N(e^{i\lambda})+A_{N+M_2+}(e^{i \lambda})- h_2(e^{i\lambda}) \bot \, L_2^{S_2}(f),\\
&\; A_{-M_1}(e^{i\lambda})+A_N(e^{i\lambda})+A_{N+M_2+}(e^{i \lambda})- h_3(e^{i\lambda}) \bot \, L_2^{S_3}(f).
\end{split} \end{equation*}

It follows from the second condition that for any $\eta \in H^{S_k}(\xi)$ the following equations should be satisfied
\begin{equation*}
\left( A_{S_k}\xi-\hat{A}_{S_k}\xi,\eta\right)=E\left[(A_{S_k}\xi-\hat{A}_{S_k}\xi)\overline{\eta}\right] =0.
 \end{equation*}

The last relations are equivalent to equations
\begin{equation*}
E\left[(A_{S_k}\xi-\hat{A}_{S_k}\xi)\overline{\xi(j)}\right] =0, \,\,\, j \in S_k.
 \end{equation*}

By using representations (3), (4), (5) and definition of the inner product in $H$ we get
\begin{equation*}
\int\limits_{-\pi}^{\pi} \left(A_{-M_1}(e^{i\lambda})+A_N(e^{i\lambda})- h_1(e^{i\lambda})\right)f(\lambda)e^{-ij\lambda }d\lambda=0,\,\,j \in S_1,
\end{equation*}
\begin{equation*}
\int\limits_{-\pi}^{\pi} \left(A_N(e^{i\lambda})+A_{N+M_2+}(e^{i\lambda})- h_2(e^{i\lambda})\right)f(\lambda)e^{-ij\lambda }d\lambda=0,\,\,j \in S_2,
\end{equation*}
\begin{equation*}
\int\limits_{-\pi}^{\pi} \left(A_{-M_1}(e^{i\lambda})+A_N(e^{i\lambda})+A_{N+M_2+}(e^{i\lambda})- h_3(e^{i\lambda})\right)f(\lambda)e^{-ij\lambda }d\lambda=0,\,\,j \in S_3.
\end{equation*}
It follows from this condition that the functions
\[(A_{-M_1}(e^{i\lambda})+A_N(e^{i\lambda})- h_1(e^{i\lambda}))f(\lambda), \; \left(A_N(e^{i\lambda})+A_{N+M_2+}(e^{i\lambda})- h_2(e^{i\lambda})\right)f(\lambda),\]
\[ \left(A_{-M_1}(e^{i\lambda})+A_N(e^{i\lambda})+A_{N+M_2+}(e^{i\lambda})- h_3(e^{i\lambda})\right)f(\lambda)\]
are of the form
\[(A_{-M_1}(e^{i\lambda})+A_N(e^{i\lambda})- h_1(e^{i\lambda}))f(\lambda)=C_N(e^{i\lambda})+C_{-M_1}(e^{i\lambda}),\]
\[\left(A_N(e^{i\lambda})+A_{N+M_2+}(e^{i\lambda})- h_2(e^{i\lambda})\right)f(\lambda)=C_N(e^{i\lambda})+C_{N+M_2+}(e^{i\lambda}),\]
\[ \left(A_{-M_1}(e^{i\lambda})+A_N(e^{i\lambda})+A_{N+M_2+}(e^{i\lambda})- h_3(e^{i\lambda})\right)f(\lambda)=C_{-M_1}(e^{i\lambda})+C_N(e^{i\lambda})+C_{N+M_2+}(e^{i\lambda})\]
where
\[C_N(e^{i\lambda})=\sum\limits_{j=0}^{N}c(j)e^{ij\lambda},\; C_{-M_1}(e^{i\lambda})=\sum\limits_{j=-\infty}^{-M_1-1}c(j)e^{ij\lambda},\; C_{N+M_2+}(e^{i\lambda})=\sum\limits_{j=N+M_2+1}^{\infty}c(j)e^{ij\lambda},
\]
and $c(j)$ are unknown coefficients that we have to find.

From the last relations we deduce that the spectral characteristics $h_k(e^{i\lambda})$ of the optimal linear estimates of the functionals $A_{S_k}\xi$ are of the form
\begin{equation} \label {5}
 h_1(e^{i\lambda})=A_{-M_1}(e^{i\lambda})+A_N(e^{i\lambda})-\frac {C_N(e^{i\lambda})+C_{-M_1}(e^{i\lambda})} {f(\lambda)},
\end{equation}
\begin{equation} \label {5_1}
 h_2(e^{i\lambda})=A_{N+M_2+}(e^{i\lambda})+A_N(e^{i\lambda})-\frac {C_N(e^{i\lambda})+C_{N+M_2+}(e^{i\lambda})} {f(\lambda)}.
\end{equation}
\begin{equation} \label {5_2}
 h_3(e^{i\lambda})=A_{-M_1}(e^{i\lambda})+A_N(e^{i\lambda})+A_{N+M_2+}(e^{i\lambda})-\frac {C_{-M_1}(e^{i\lambda})+C_N(e^{i\lambda})+C_{N+M_2+}(e^{i\lambda})} {f(\lambda)}.
\end{equation}

To find equations for calculation of the unknown coefficients $c(j)$ we use the decomposition of the function $(f(\lambda))^{-1}$ into the Fourier series
\begin{equation} \label{decomp}
\frac {1} {f(\lambda)}=\sum\limits_{m=-\infty}^{\infty}b(m) e^{im\lambda},
\end{equation}
where $b(m)$ are the Fourier coefficients of the function $(f(\lambda))^{-1}$.

Inserting (9) into (6), (7),  (8) we obtain the following representations of the spectral characteristics
\begin{equation} \label{spchar}
 h_1(e^{i\lambda})=\sum\limits_{j=-\infty}^{-M_1-1}a(j)e^{ij\lambda}+\sum\limits_{j=0}^{N}a(j)e^{ij\lambda }-\left(\sum\limits_{m=-\infty}^{\infty}b(m) e^{im\lambda}\right)\left(\sum\limits_{j=0}^{N}c(j)e^{ij\lambda}+\sum\limits_{j=-\infty}^{-M_1-1}c(j)e^{ij\lambda}\right),
\end{equation}
\[ h_2(e^{i\lambda})=\sum\limits_{j=N+M_2+1}^{\infty}a(j)e^{ij\lambda}+\sum\limits_{j=0}^{N}a(j)e^{ij\lambda }-\]
\begin{equation} \label{spchar_1}
-\left(\sum\limits_{m=-\infty}^{\infty}b(m) e^{im\lambda}\right)\left(\sum\limits_{j=0}^{N}c(j)e^{ij\lambda}+\sum\limits_{j=N+M_2+1}^{\infty}c(j)e^{ij\lambda}\right),
\end{equation}
\[ h_3(e^{i\lambda})=\sum\limits_{j=-\infty}^{-M_1-1}a(j)e^{ij\lambda}+\sum\limits_{j=0}^{N}a(j)e^{ij\lambda }+\sum\limits_{j=N+M_2+1}^{\infty}a(j)e^{ij\lambda}-\]
\begin{equation} \label{spchar_2}
-\left(\sum\limits_{m=-\infty}^{\infty}b(m) e^{im\lambda}\right)\left(\sum\limits_{j=-\infty}^{-M_1-1}c(j)e^{ij\lambda}+\sum\limits_{j=0}^{N}c(j)e^{ij\lambda}+\sum\limits_{j=N+M_2+1}^{\infty}c(j)e^{ij\lambda}\right),
\end{equation}

It follows from the first conditions $h_k(e^{i\lambda}) \in L_2^{S_k}(F)$ that the Fourier coefficients of the functions  $h_k(e^{i\lambda})$ are equal to zero for $j \in Z\backslash S_k$ respectively, that is,
\[\int\limits_{-\pi}^{\pi} h_k(e^{i\lambda})e^{-ij\lambda }d\lambda=0,\, j \in Z\backslash S_k.
\]

Using the last relations and (10), (11), (12) we get the following system of equations that determine the unknown coefficients $c(j),j \in Z\backslash S_1$,
\begin{equation} \label{koef_1}
\begin{split}
a(0)=&\sum\limits_{j=0}^{N}c(j)b(-j)+\sum\limits_{j=M_1+1}^{\infty}c(-j)b(j);\\
a(1)=&\sum\limits_{j=0}^{N}c(j)b(-j+1)+\sum\limits_{j=M_1+1}^{\infty}c(-j)b(j+1);\\
\ldots\\
a(N)=&\sum\limits_{j=0}^{N}c(j)b(-j+N)+\sum\limits_{j=M_1+1}^{\infty}c(-j)b(j+N);\\
a(-M_1-1)=\sum\limits_{j=0}^{N}&c(j)b(-M_1-1-j)+\sum\limits_{j=M_1+1}^{\infty}c(-j)b(j-M_1-1);\\
a(-M_1-2)=\sum\limits_{j=0}^{N}&c(j)b(-M_1-2-j)+\sum\limits_{j=M_1+1}^{\infty}c(-j)b(j-M_1-2);\\
\ldots\\
\end{split} \end{equation}
the following system of equations that determine the unknown coefficients $c(j),j \in Z\backslash S_2$
\begin{equation} \label{koef_2}
\begin{split}
a(0)=&\sum\limits_{j=0}^{N}c(j)b(-j)+\sum\limits_{j=N+M_2+1}^{\infty}c(j)b(-j);\\
a(1)=&\sum\limits_{j=0}^{N}c(j)b(-j+1)+\sum\limits_{j=N+M_2+1}^{\infty}c(j)b(-j+1);\\
\ldots\\
a(N)=&\sum\limits_{j=0}^{N}c(j)b(-j+N)+\sum\limits_{j=N+M_2+1}^{\infty}c(j)b(-j+N);\\
a(N+M_2+1)=\sum\limits_{j=0}^{N}&c(j)b(N+M_2+1-j)+\sum\limits_{j=N+M_2+1}^{\infty}c(j)b(N+M_2+1-j);\\
a(N+M_2+2)=\sum\limits_{j=0}^{N}&c(j)b(N+M_2+2-j)+\sum\limits_{j=N+M_2+1}^{\infty}c(j)b(N+M_2+2-j);\\
\ldots\\
\end{split} \end{equation}
and the following system of equations that determine the unknown coefficients $c(j),j \in Z\backslash S_3$
\[\begin{split}
\ldots\\
a(-M_1-2)=\sum\limits_{j=M_1+1}^{\infty}c(-j)b(j-M_1-2)+\sum\limits_{j=0}^{N}c(j)b(-M_1-2-j)+\sum\limits_{j=N+M_2+1}^{\infty}c(j)b(-M_1-2-j);\\
a(-M_1-1)=\sum\limits_{j=M_1+1}^{\infty}c(-j)b(j-M_1-1)+\sum\limits_{j=0}^{N}c(j)b(-M_1-1-j)+\sum\limits_{j=N+M_2+1}^{\infty}c(j)b(-M_1-1-j);\\
\end{split}\]
\begin{equation} \label{koef_3}
\begin{split}
a(0)=\sum\limits_{j=M_1+1}^{\infty}c(-j)b(j)+\sum\limits_{j=0}^{N}c(j)b(-j)+\sum\limits_{j=N+M_2+1}^{\infty}c(j)b(-j);\\
a(1)=\sum\limits_{j=M_1+1}^{\infty}c(-j)b(j+1)+\sum\limits_{j=0}^{N}c(j)b(-j+1)+\sum\limits_{j=N+M_2+1}^{\infty}c(j)b(-j+1);\\
\ldots\\
a(N)=\sum\limits_{j=M_1+1}^{\infty}c(-j)b(j+N)+\sum\limits_{j=0}^{N}c(j)b(-j+N)+\sum\limits_{j=N+M_2+1}^{\infty}c(j)b(-j+N);\\
\end{split} \end{equation}
\[\begin{split}
a(N+M_2+1)=\sum\limits_{j=M_1+1}^{\infty}c(-j)b(N+M_2+1+j)+\sum\limits_{j=0}^{N}c(j)b(N+M_2+1-j)+\\
+\sum\limits_{j=N+M_2+1}^{\infty}c(j)b(N+M_2+1-j);\\
a(N+M_2+2)=\sum\limits_{j=M_1+1}^{\infty}c(-j)b(N+M_2+2+j)+\sum\limits_{j=0}^{N}c(j)b(N+M_2+2-j)+\\
+\sum\limits_{j=N+M_2+1}^{\infty}c(j)b(N+M_2+2-j);\\
\ldots\\
\end{split}\]

Denote by
\[\vec{\bold{a}}_{1}^\top=(\vec{a}_N^\top,\vec{a}_{-M_1-1}^\top),\; \vec{a}_N^\top=(a(0),a(1),\ldots,a(N)),\; \vec{a}_{-M_1-1}^\top=(a(-M_1-1),a(-M_1-2),\ldots),\]
\[\vec{\bold{a}}_{2}^\top=(\vec{a}_N^\top,\vec{a}_{N+M_2+1}^\top),\; \vec{a}_{N+M_2+1}^\top=(a(N+M_2+1),a(N+M_2+2),\ldots).\]

Let $B_{S_1}$ be a matrix
\begin {equation*}
B_{S_1}= \left( \begin{array}{cc}
B_{1}& B_{2}\\
B_{3}& B_{4}
\end{array}\right),
\end{equation*}
where $B_{1}$, $B_{2}$, $B_{3}$, $B_{4}$ are matrices constructed from the Fourier coefficients of the function $(f(\lambda))^{-1}$:
\[B_1(i,j)=b(i-j),\; 0\le i,j\le N,\]
\[B_2(i,j)=b(M_1+1+i+j),\; 0\le i\le N,\;0\le j<\infty,\]
\[B_3(i,j)=b(-M_1-1-i-j),\; 0\le j\le N,\;0\le i<\infty,\]
\[B_4(i,j)=b(j-i),\; 0\le i,j<\infty.\]

Let $B_{S_2}$ be a matrix
\begin {equation*}
B_{S_2}= \left( \begin{array}{cc}
B_1& B_5\\
B_6& B_7
\end{array}\right),
\end{equation*}
where $B_5$, $B_6$, $B_7$ are matrices constructed from the Fourier coefficients of the function $(f(\lambda))^{-1}$:
\[B_5(i,j)=b(-N-M_2-1+i-j),\; 0\le i\le N,\;0\le j<\infty,\]
\[B_6(i,j)=b(N+M_2+1+i-j),\; 0\le j\le N,\;0\le i<\infty,\]
\[B_7(i,j)=b(i-j),\; 0\le i,j<\infty.\]
Making use the introduced notations we can write formulas (13}), (14) in the form of equations
\begin{equation} \label{rivn}
\vec{\bold{a}}_{1}=B_{S_1}\vec{\bold{c}}_{1},\; \vec{\bold{a}}_{2}=B_{S_2}\vec{\bold{c}}_{2},
\end{equation}
and  formulas (15) in the form of equations
\begin{equation} \label{rivn1}
\vec{a}_N=B_1\vec{c}_N+B_2\vec{c}_{-M_1-1}+B_5\vec{c}_{N+M_2+1},\;
\vec{a}_{-M_1-1}=B_3\vec{c}_N+B_4\vec{c}_{-M_1-1}+B_8\vec{c}_{N+M_2+1},
\end{equation}
\begin{equation} \label{rivn2}
\vec{a}_{N+M_2+1}=B_6\vec{c}_N+B_9\vec{c}_{-M_1-1}+B_7\vec{c}_{N+M_2+1},
\end{equation}
where $\vec{\bold{c}}_{1}^\top=(\vec{c}_N^\top,\vec{c}_{-M_1-1}^\top)$, $\vec{\bold{c}}_{2}^\top=(\vec{c}_N^\top,\vec{c}_{N+M_2+1}^\top)$ are vectors constructed from the unknown coefficients $c(j)$:
\[\vec{c}_N^\top=(c(0),c(1),\ldots,c(N)),\; \vec{c}_{-M_1-1}^\top=(c(-M_1-1),c(-M_1-2),\ldots),\;\] \[\vec{c}_{N+M_2+1}^\top=(c(N+M_2+1),c(N+M_2+2),\ldots),\]
and $B_8$, $B_9$ are matrices constructed from the Fourier coefficients of the function $(f(\lambda))^{-1}$:
\[B_8(i,j)=b(-N-M_2-1-j-M_1-1-i),\; 0\le i,j<\infty,\]
\[B_9(i,j)=b(N+M_2+1+i+M_1+1+j),\; 0\le i,j<\infty.\]

Since the matrix $B_{S_k}$ is reversible [16, 21], we get the formulas
\begin{equation} \label{rivn2}
\vec{\bold{c}}_{1}=B_{S_1}^{-1}\vec{\bold{a}}_{1},\; \vec{\bold{c}}_{2}=B_{S_2}^{-1}\vec{\bold{a}}_{2}.
\end{equation}

 Hence, the unknown coefficients $c(j), j \in Z\backslash S_1$, and $c(j), j \in Z\backslash S_2$, are calculated  by the formulas
\[c(j)=\left(B_{S_1}^{-1}\vec{\bold{a}}_{1}\right)(j),\; j \in Z\backslash S_1,\]
\[c(j)=\left(B_{S_2}^{-1}\vec{\bold{a}}_{2}\right)(j),\; j \in Z\backslash S_2.\]
where $\left(B_{S_1}^{-1}\vec{\bold{a}}_{1}\right)(j)$ and $\left(B_{S_2}^{-1}\vec{\bold{a}}_{2}\right)(j)$ are the  $j$ components of the vectors $B_{S_1}^{-1}\vec{\bold{a}}_{1}$ and $B_{S_2}^{-1}\vec{\bold{a}}_{2}$ respectively, and the formulas for calculating  the spectral characteristics of the estimates $\hat{A}_{S_1}\xi$, $\hat{A}_{S_2}\xi$ are of the form
\[h_1(e^{i\lambda})=A_{-M_1}(e^{i\lambda})+A_N(e^{i\lambda})-\]
\begin{equation} \label{ch_1}
-\left(\sum\limits_{m=-\infty}^{\infty}b(m) e^{im\lambda}\right)\left(\sum\limits_{j=0}^{N}\left(B_{S_1}^{-1}\vec{\bold{a}}_{1}\right)(j)e^{ij\lambda}+\sum\limits_{j=-\infty}^{-M_1-1}\left(B_{S_1}^{-1}\vec{\bold{a}}_{1}\right)(j)e^{ij\lambda}\right),
\end{equation}
\[h_2(e^{i\lambda})=A_{N+M_2+}(e^{i\lambda})+A_N(e^{i\lambda})-\]
\begin{equation} \label{ch_2}
-\left(\sum\limits_{m=-\infty}^{\infty}b(m) e^{im\lambda}\right)\left(\sum\limits_{j=0}^{N}\left(B_{S_2}^{-1}\vec{\bold{a}}_{2}\right)(j)e^{ij\lambda}+\sum\limits_{j=N+M_2+1}^{\infty}\left(B_{S_2}^{-1}\vec{\bold{a}}_{2}\right)(j)e^{ij\lambda}\right).
\end{equation}

The mean square errors of the estimates of the functionals can be calculated by the formulas
\[\Delta(h_1;f)=\frac{1}{2 \pi} \int_{- \pi}^{ \pi} \, \left|C_{N} (e^{i \lambda} ))+C_{-M_1} (e^{i \lambda})\right|^{2} f^{-1}( \lambda ) d \lambda=\]
\[=\frac{1}{2 \pi}\int\limits_{-\pi}^{\pi}\left(\sum\limits_{j=0}^{N}c(j)e^{ij\lambda}+\sum\limits_{j=-\infty}^{-M_1-1}c(j)e^{ij\lambda}\right)\left(\sum\limits_{j=0}^{N}\overline{c(j)}e^{-ij\lambda}+\sum\limits_{j=-\infty}^{-M_1-1}\overline{c(j)}e^{-ij\lambda}\right)\left(\sum\limits_{m=-\infty}^{\infty}b(m)e^{im\lambda}\right)d\lambda=\]
\begin {equation} \label {6}
=\left<\vec{\bold{c}}_{1},B_{S_1}\vec{\bold{c}}_{1}\right>=\left<B_{S_1}^{-1}\vec{\bold{a}}_{1}, \vec{\bold{a}}_{1}\right>,
\end{equation}
\[\Delta(h_2;f)=\frac{1}{2 \pi} \int_{- \pi}^{ \pi} \, \left|C_{N} (e^{i \lambda} ))+C_{N+M_2+} (e^{i \lambda})\right|^{2} f^{-1}( \lambda ) d \lambda=\]
\[=\frac{1}{2 \pi}\int\limits_{-\pi}^{\pi}\left(\sum\limits_{j=0}^{N}c(j)e^{ij\lambda}+\sum\limits_{j=N+M_2+1}^{\infty}c(j)e^{ij\lambda}\right)\left(\sum\limits_{j=0}^{N}\overline{c(j)}e^{-ij\lambda}+\sum\limits_{j=N+M_2+1}^{\infty}\overline{c(j)}e^{-ij\lambda}\right)\left(\sum\limits_{m=-\infty}^{\infty}b(m)e^{im\lambda}\right)d\lambda=\]
\begin {equation} \label {6_1}
=\left<\vec{\bold{c}}_{2},B_{S_2}\vec{\bold{c}}_{2}\right>=\left<B_{S_2}^{-1}\vec{\bold{a}}_{2}, \vec{\bold{a}}_{2}\right>,
\end{equation}
\[\Delta(h_3;f)=\frac{1}{2 \pi} \int_{- \pi}^{ \pi} \, \left|C_{-M_1} (e^{i \lambda})+C_{N} (e^{i \lambda} ))+C_{N+M_2+} (e^{i \lambda})\right|^{2} f^{-1}( \lambda ) d \lambda=\]
\[=\frac{1}{2 \pi}\int\limits_{-\pi}^{\pi}\left(\sum\limits_{j=-\infty}^{-M_1-1}c(j)e^{ij\lambda}+\sum\limits_{j=0}^{N}c(j)e^{ij\lambda}+\sum\limits_{j=N+M_2+1}^{\infty}c(j)e^{ij\lambda}\right)\times\]
\[\times\left(+\sum\limits_{j=-\infty}^{-M_1-1}\overline{c(j)}e^{-ij\lambda}+\sum\limits_{j=0}^{N}\overline{c(j)}e^{-ij\lambda}+\sum\limits_{j=N+M_2+1}^{\infty}\overline{c(j)}e^{-ij\lambda}\right)\left(\sum\limits_{m=-\infty}^{\infty}b(m)e^{im\lambda}\right)d\lambda=\]
\begin {equation} \label {6_2}
=\left<\vec{c}_{N},\vec{a}_{N}\right>+\left<\vec{c}_{-M_1-1},\vec{a}_{-M_1-1}\right>+\left<\vec{c}_{N+M_2+1},\vec{a}_{N+M_2+1}\right>,
\end{equation}
where $\left<\cdot,\cdot\right> $\,is the inner product.

Let us summarize our results and present them in the form of a theorem.

\begin{thm}
Let $\xi(j)$ be a stationary stochastic sequence with the spectral density $f(\lambda)$ that satisfies the minimality condition (1). The mean square errors  $\Delta(h_k,f)$ and the spectral characteristics $h_k(e^{i\lambda})$ of the optimal linear estimates $\hat{A}_{S_k}\xi$ of the functionals $A_{S_k}\xi$ from observations of the sequence $\xi(j)$ at points $j\in S_k$ respectively, where
\[S_1=\{-M_1, \ldots, -1\}\cup \{N+1, N+2, \ldots\}, \; S_2=\{\ldots, -2, -1\}\cup  \{N+1, \ldots, N+M_2\},\]
\[S_3=\{-M_1, \ldots, -1\}\cup \{N+1, \ldots, N+M_2\},\]
can be calculated by  formulas (22) -- (24), (20), (21), (8).
\end{thm}

Consider the problem of the mean-square optimal estimation of the linear functionals
\[A_{S_4}\xi=\sum\limits_{j=-\-M_1-N_1}^{-M_1-1}a(j)\xi(j)+\sum\limits_{j=0}^{N}a(j)\xi(j),\; A_{S_5}\xi=\sum\limits_{j=0}^{N}a(j)\xi(j)+\sum\limits_{j=N+M_2+1}^{N+M_2+N_2}a(j)\xi(j),\]
\[A_{S_6}\xi=\sum\limits_{j=-\-M_1-N_1}^{-M_1-1}a(j)\xi(j)+
\sum\limits_{j=0}^{N}a(j)\xi(j)+\sum\limits_{j=N+M_2+1}^{N+M_2+N_2}a(j)\xi(j),\]
which depend on the unknown values of a stochastic stationary sequence $\xi(j),j\in\mathbb{Z},$  from observations of the sequence at points $j\in S_k$ respectively, where
\[S_4=Z\backslash(\{-M_1-N_1, \ldots, -M_1-1\}\cup\{0, \ldots, N\}),\]
\[S_5=Z\backslash(\{0, \ldots, N\}\cup\{N+M_2+1, \ldots, N+M_2+N_1\}),\]
 \[S_6=Z\backslash(\{-M_1-N_1, \ldots, -M_1-1\}\cup\{0, \ldots, N\}\cup\{N+M_2+1, \ldots, N+M_2+N_1\}).\]

Denote by
\[\vec{\bold{a}}_{4}^\top=(\vec{a}_N^\top,\vec{a}_{-M_1-1}^{N_1\top}),\; \vec{\bold{a}}_{5}^\top=(\vec{a}_N^\top,\vec{a}_{N+M_2+1}^{N_2\top}),\; \vec{\bold{a}}_{6}^\top=(\vec{a}_N^\top, \vec{a}_{-M_1-1}^{N_1\top},\vec{a}_{N+M_2+1}^{N_2\top}),\]
\[\vec{a}_N^\top=(a(0),a(1),\ldots,a(N)),\; \vec{a}_{-M_1-1}^{N_1\top}=(a(-M_1-1),a(-M_1-2),\ldots,a(-M_1-N_1)),\]
\[\vec{a}_{N+M_2+1}^{N_2\top}=(a(N+M_2+1),a(N+M_2+2),\ldots,a(N+M_2+N_2)).\]

Let $B_{S_4}$ be a matrix
\begin {equation*}
B_{S_4}= \left( \begin{array}{cc}
B_{1}& B_{2}\\
B_{3}& B_{4}
\end{array}\right),
\end{equation*}
where $B_{1}$, $B_{2}$, $B_{3}$, $B_{4}$ are matrices constructed from the Fourier coefficients of the function $(f(\lambda))^{-1}$:
\[B_1(i,j)=b(i-j),\; 0\le i,j\le N,\]
\[B_2(i,j)=b(M_1+1+i+j),\; 0\le i\le N,\;0\le j<N_1,\]
\[B_3(i,j)=b(-M_1-1-i-j),\; 0\le j\le N,\;0\le i<N_1,\]
\[B_4(i,j)=b(j-i),\; 0\le i,j<N_1.\]

Let $B_{S_5}$ be a matrix
\begin {equation*}
B_{S_5}= \left( \begin{array}{cc}
B_1& B_5\\
B_6& B_7
\end{array}\right),
\end{equation*}
where $B_5$, $B_6$, $B_7$ are matrices constructed from the Fourier coefficients of the function $(f(\lambda))^{-1}$:
\[B_5(i,j)=b(-N-M_2-1+i-j),\; 0\le i\le N,\;0\le j<N_2,\]
\[B_6(i,j)=b(N+M_2+1+i-j),\; 0\le j\le N,\;0\le i<N_2,\]
\[B_7(i,j)=b(i-j),\; 0\le i,j<N_2.\]

Let $B_{S_6}$ be a matrix
\begin {equation*}
B_{S_5}= \left( \begin{array}{ccc}
B_1&B_2& B_5\\
B_{3}& B_{4}&B_8\\
B_6&B_9& B_7
\end{array}\right),
\end{equation*}
where $B_8$, $B_9$ are matrices constructed from the Fourier coefficients of the function $(f(\lambda))^{-1}$:
\[B_8(i,j)=b(-N-M_2-1-j-M_1-1-i),\; 0\le i<N_1,\; 0\le j<N_2\]
\[B_9(i,j)=b(N+M_2+1+i+M_1+1+j),\; 0\le i<N_1,\; 0\le j<N_2.\]

\begin{nas}
Let $\xi(j)$ be a stationary stochastic sequence with the spectral density $f(\lambda)$ that satisfies the minimality condition (1). The mean square errors  $\Delta(h_k,f)$ and the spectral characteristics $h_k(e^{i\lambda})$ of the optimal linear estimates $\hat{A}_{S_k}\xi$ of the functionals $A_{S_k}\xi$ from observations of the sequence $\xi(j)$ at points $j\in S_k$ respectively, where k=4, 5, 6, can be calculated by the following formulas
\[h_4(e^{i\lambda})=\sum\limits_{j=-M_1-N_1}^{-M_1-1}a(j)e^{ij\lambda}+\sum\limits_{j=0}^{N}a(j)e^{ij\lambda }-\]
\begin{equation} \label{ch_11}
-\left(\sum\limits_{m=-\infty}^{\infty}b(m) e^{im\lambda}\right)\left(\sum\limits_{j=0}^{N}\left(B_{S_4}^{-1}\vec{\bold{a}}_{4}\right)(j)e^{ij\lambda}+\sum\limits_{j=-M_1-N_1}^{-M_1-1}\left(B_{S_4}^{-1}\vec{\bold{a}}_{4}\right)(j)e^{ij\lambda}\right),
\end{equation}
\begin {equation} \label {7}
\Delta(h_4;f)=\left<B_{S_4}^{-1}\vec{\bold{a}}_{4}, \vec{\bold{a}}_{4}\right>,
\end{equation}
\[h_5(e^{i\lambda})=\sum\limits_{j=0}^{N}a(j)e^{ij\lambda }+\sum\limits_{j=N+M_2+1}^{N+M_2+N_2}a(j)e^{ij\lambda}-\]
\begin{equation} \label{ch_12}
-\left(\sum\limits_{m=-\infty}^{\infty}b(m) e^{im\lambda}\right)\left(\sum\limits_{j=0}^{N}\left(B_{S_5}^{-1}\vec{\bold{a}}_{5}\right)(j)e^{ij\lambda}+\sum\limits_{j=N+M_2+1}^{N+M_2+N_2}\left(B_{S_5}^{-1}\vec{\bold{a}}_{5}\right)(j)e^{ij\lambda}\right).
\end{equation}
\begin {equation} \label {7_1}
\Delta(h_5;f)=\left<B_{S_5}^{-1}\vec{\bold{a}}_{5}, \vec{\bold{a}}_{5}\right>,
\end{equation}
\[ h_6(e^{i\lambda})=\sum\limits_{j=-M_1-N_1}^{-M_1-1}a(j)e^{ij\lambda}+\sum\limits_{j=0}^{N}a(j)e^{ij\lambda }+\sum\limits_{j=N+M_2+1}^{N+M_2+N_2}a(j)e^{ij\lambda}-\left(\sum\limits_{m=-\infty}^{\infty}b(m)e^{im\lambda}\right)\times\]
\begin{equation} \label{ch_22}
\times\left(\sum\limits_{j=0}^{N}\left(B_{S_6}^{-1}\vec{\bold{a}}_{6}\right)(j)e^{ij\lambda}+\sum\limits_{j=-M_1-N-1}^{-M_1-1}\left(B_{S_6}^{-1}\vec{\bold{a}}_{6}\right)(j)e^{ij\lambda}+\sum\limits_{j=N+M_2+1}^{N+M_2+N_2}\left(B_{S_6}^{-1}\vec{\bold{a}}_{6}\right)(j)e^{ij\lambda}\right),
\end{equation}
\begin {equation} \label {7_2}
\Delta(h_6;f)=\left<B_{S_6}^{-1}\vec{\bold{a}}_{6}, \vec{\bold{a}}_{6}\right>.
\end{equation}
\end{nas}

\begin{exm}
Let's investigate the problem of linear estimation of the functional
\[A_{S_4}\xi=a(-5)\xi(-5)+a(-4)\xi(-4)+a(-3)\xi(-3)+a(0)\xi(0)+a(1)\xi(1)\]
based on observations of the sequence $\xi(j)$ at points $j\in Z\backslash \{-5, -4, -3, 0, 1\}$. Consider this problem for the stationary sequence wiyh the spectral density
\[f(\lambda)=\frac{1}{\left|1-\alpha e^{-i\lambda}\right|^2},\; |\alpha|<1.\]
In this case the decomposition of the function $f^{-1}(\lambda)$ into the Fourier series is as follows
\[f^{-1}(\lambda)=\left|1-\alpha e^{-i\lambda}\right|^2=1+|\alpha|^2-\overline{\alpha}e^{i\lambda}-\alpha e^{-i\lambda}.\]
The coefficients which determine the spectral characteristic and the mean-square error of the functional are the following
\[b(0)=1+|\alpha|^2,\; b(1)=-\overline{\alpha},\; b(-1)=-\alpha.\]
The spectral characteristic  of the optimal estimate is calculated by the formula
\[h_4(e^{i\lambda})=a(-5)e^{-5i\lambda}+a(-4)e^{-4i\lambda}+a(-3)e^{-3i\lambda}+a(0)+a(1)e^{i\lambda}-\]
\[-\left(c(-5)e^{-5i\lambda}+c(-4)e^{-4i\lambda}+c(-3)e^{-3i\lambda}+c(0)+c(1)e^{i\lambda}\right)\left(1+|\alpha|^2-\overline{\alpha}e^{i\lambda}-\alpha e^{-i\lambda}\right)=\]
\[=\overline{\alpha}c(1)e^{2i\lambda}+\alpha c(0)e^{-i\lambda}+\overline{\alpha}c(-3)e^{-2i\lambda}+\alpha c(-5)e^{-6i\lambda},\]
where the unknown coefficients $c(k)$ are determined by formulas
\[c(0)=\frac{a(0)(1+|\alpha|^2)+\alpha a(1)}{1+|\alpha|^2+|\alpha|^4},\; c(1)=\frac{\overline{\alpha}a(0)+a(1)(1+|\alpha|^2)}{1+|\alpha|^2+|\alpha|^4},\]
\[c(-3)=\frac{(1+|\alpha|^2+|\alpha|^4)a(-3)+\overline{\alpha}(1+|\alpha|^2)a(-4)+\overline{\alpha}^2a(-5)}{(1+|\alpha|^2)(1+|\alpha|^4)},\]
\[c(-4)=\frac{\alpha a(-3)+(1+|\alpha|^2)a(-4)+\overline{\alpha}a(-5)}{1+|\alpha|^4},\]
\[c(-5)=\frac{\alpha^2 a(-3)+\alpha(1+|\alpha|^2)a(-4)+(1+|\alpha|^2+|\alpha|^4)a(-5)}{(1+|\alpha|^2)(1+|\alpha|^4)}.\]
The mean-square error is of the form
\[\Delta_4(f)=((1+|\alpha|^2)(a^2(0)+a^2(1))+(\alpha+\overline{\alpha})a(0)a(1))(1+|\alpha|^2+|\alpha|^4)^{-1}+\]
\[+((1+|\alpha|^2+|\alpha|^4)(a^2(-3)+a^2(-4)+a^2(-5))+|\alpha|^2a^2(-4)+(\alpha^2+\overline{\alpha}^2)a(-3)a(-5)+\]
\[+(\alpha+\overline{\alpha})(1+|\alpha|^2)(a(-3)a(-4)+a(-4)a(-5)))((1+|\alpha|^2)(1+|\alpha|^4))^{-1}.\]
Consider the problem of linear estimation of the functional
\[A_{S_5}\xi=a(0)\xi(0)+a(1)\xi(1)+a(4)\xi(4)+a(5)\xi(5)+a(6)\xi(6)\]
based on observations of the sequence $\xi(j)$ at points $j\in Z\backslash \{0, 1, 4, 5, 6\}$. The spectral characteristic  of the optimal estimate is calculated by the formula
\[h_5(e^{i\lambda})=a(0)+a(1)e^{i\lambda}+a(4)e^{4i\lambda}+a(5)e^{5i\lambda}+a(6)e^{6i\lambda}-\]
\[-\left(c(0)+c(1)e^{i\lambda}+c(4)e^{4i\lambda}+c(5)e^{5i\lambda}+c(6)e^{6i\lambda}\right)\left(1+|\alpha|^2-\overline{\alpha}e^{i\lambda}-\alpha e^{-i\lambda}\right)=\]
\[=\alpha c(0)e^{-i\lambda}+\overline{\alpha}c(1)e^{2i\lambda}+\alpha c(4)e^{3i\lambda}+\overline{\alpha} c(6)e^{7i\lambda},\]
where the unknown coefficients $c(k)$ are determined by formulas
\[c(0)=\frac{a(0)(1+|\alpha|^2)+\alpha a(1)}{1+|\alpha|^2+|\alpha|^4},\; c(1)=\frac{\overline{\alpha}a(0)+a(1)(1+|\alpha|^2)}{1+|\alpha|^2+|\alpha|^4},\]
\[c(4)=\frac{(1+|\alpha|^2+|\alpha|^4)a(4)+\alpha(1+|\alpha|^2)a(5)+\alpha^2a(6)}{(1+|\alpha|^2)(1+|\alpha|^4)},\]
\[c(5)=\frac{\overline{\alpha} a(4)+(1+|\alpha|^2)a(5)+\alpha a(6)}{1+|\alpha|^4},\]
\[c(6)=\frac{\overline{\alpha}^2 a(4)+\overline{\alpha}(1+|\alpha|^2)a(5)+(1+|\alpha|^2+|\alpha|^4)a(6)}{(1+|\alpha|^2)(1+|\alpha|^4)}.\]
The mean-square error is of the form
\[\Delta_5(f)=((1+|\alpha|^2)(a^2(0)+a^2(1))+(\alpha+\overline{\alpha})a(0)a(1))(1+|\alpha|^2+|\alpha|^4)^{-1}+\]
\[+((1+|\alpha|^2+|\alpha|^4)(a^2(4)+a^2(5)+a^2(6))+|\alpha|^2a^2(5)+(\alpha^2+\overline{\alpha}^2)a(4)a(6)+\]
\[+(\alpha+\overline{\alpha})(1+|\alpha|^2)(a(5)a(6)+a(4)a(5)))((1+|\alpha|^2)(1+|\alpha|^4))^{-1}.\]
Consider the problem of linear estimation of the functional
\[A_{S_6}\xi=a(-5)\xi(-5)+a(-4)\xi(-4)+a(-3)\xi(-3)+a(0)\xi(0)+a(1)\xi(1)+a(4)\xi(4)+a(5)\xi(5)+a(6)\xi(6)\]
based on observations of the sequence $\xi(j)$ at points $j\in Z\backslash \{-5, -4, -3, 0, 1, 4, 5, 6\}$. The spectral characteristic  of the optimal estimate is calculated by the formula
\[h_6(e^{i\lambda})=\overline{\alpha}c(-3)e^{-2i\lambda}+\alpha c(-5)e^{-6i\lambda}+\alpha c(0)e^{-i\lambda}+\overline{\alpha}c(1)e^{2i\lambda}+\alpha c(4)e^{3i\lambda}+\overline{\alpha} c(6)e^{7i\lambda}.\]
The mean-square error is of the form
\[\Delta_6(f)=\Delta_4(f)+\Delta_5(f)-((1+|\alpha|^2)(a^2(0)+a^2(1))+(\alpha+\overline{\alpha})a(0)a(1))(1+|\alpha|^2+|\alpha|^4)^{-1}.\]

\end{exm}

\section{Minimax-robust method of interpolation}

The traditional methods of estimation of the functionals $A_{S_k}\xi$ which depends on unknown values of a stationary stochastic sequence $\xi(j)$ can be applied in the case  where the spectral density $f(\lambda)$ of the considered stochastic sequence $\xi(j)$ is exactly known. In practise, however, we do not have complete information on spectral density of the sequence. For this reason  we apply the minimax(robust) method of estimation of the functionals $A_{S_k}\xi$, that is we find an estimate that minimizes the maximum of the mean square errors for all spectral densities from the given class of admissible spectral densities $D$.

\begin{ozn}
For a given class of spectral densities $D$ spectral densities $f_k^0(\lambda)\in D$ are called the least favourable in $D$ for the optimal linear estimation of the functionals $A_{S_k}\xi$ if the following relations hold true
$$\Delta\left(f_k^0\right)=\Delta\left(h_k\left(f_k^0\right);f_k^0\right)=\max\limits_{f\in D}\Delta\left(h_k\left(f\right);f\right).$$
\end{ozn}

\begin{ozn}
For a given class of spectral densities $D$ the spectral characteristics $h_k^0(e^{i\lambda})$ of the optimal linear estimates of the functionals $A_{S_k}\xi$ are called minimax-robust if
$$h_k^0(e^{i\lambda})\in H_D^k= \bigcap\limits_{f\in D} L_2^{S_k}(f),$$
$$\min\limits_{h\in H_D^k}\max\limits_{f\in D}\Delta\left(h;f\right)=\sup\limits_{f\in D}\Delta\left(h_k^0;f\right).$$
\end{ozn}

It follows from the introduced definitions and the obtained formulas that the following statement holds true.

\begin{lem} The spectral densities $f_k^0(\lambda)\in D$ are the least favourable in the class of admissible spectral densities  $D$ for the optimal linear estimates of the functionals $A_{S_k}\xi$ if the Fourier coefficients of the functions $(f_k^0(\lambda))^{-1}$ define matrices $B_{S_k}^0$ that are solutions to the optimization problems
\begin{equation} \label{extrem}
\max\limits_{f\in D}\left<B_{S_k}^{-1}\vec{\bold{a}}_{k},\,\vec{\bold{a}}_{k}\right>=\left<B_{S_k}^{0-1}\vec{\bold{a}}_{k},\,\vec{\bold{a}}_{k}\right>.
\end{equation}
The minimax spectral characteristics $h_k^0=h_k(f_k^0)$ can be calculated by the formulas (20), (21), (25), (27), (29) if $h_k(f_k^0) \in H_D^k$.
\end{lem}

\begin{lem} The spectral densities $f_3^0(\lambda)\in D$ are the least favourable in the class of admissible spectral densities  $D$ for the optimal linear estimates of the functionals $A_{S_3}\xi$ if the Fourier coefficients of the functions $(f_3^0(\lambda))^{-1}$ define matrices $B_{k}^0$ that are solutions to the optimization problems
\[\max\limits_{f\in D}\left(\left<\vec{c}_{N},\vec{a}_{N}\right>+\left<\vec{c}_{-M_1-1},\vec{a}_{-M_1-1}\right>+\left<\vec{c}_{N+M_2+1},\vec{a}_{N+M_2+1}\right>\right)=\]
\begin{equation} \label{extrem1}
=\left<\vec{c}_{N}^0,\vec{a}_{N}\right>+\left<\vec{c}_{-M_1-1}^0,\vec{a}_{-M_1-1}\right>+\left<\vec{c}_{N+M_2+1}^0,\vec{a}_{N+M_2+1}\right>,
\end{equation}
where
\[\vec{a}_N=B_1^0\vec{c}_N^0+B_2^0\vec{c}_{-M_1-1}^0+B_5^0\vec{c}_{N+M_2+1}^0,\;
\vec{a}_{-M_1-1}=B_3^0\vec{c}_N^0+B_4^0\vec{c}_{-M_1-1}^0+B_8^0\vec{c}_{N+M_2+1}^0,\]
\[\vec{a}_{N+M_2+1}=B_6^0\vec{c}_N^0+B_9^0\vec{c}_{-M_1-1}^0+B_7^0\vec{c}_{N+M_2+1}^0.\]
The minimax spectral characteristics $h_3^0=h_3(f_3^0)$ can be calculated by the formulas (8) if $h_3(f_3^0) \in H_D^3$.
\end{lem}

The least favourable spectral densities $f_k^0$ and the minimax spectral characteristics $h_k^0$ form saddle points of the function $\Delta \left(h;f\right)$ on the sets $H_D^k\times D$. The saddle point inequalities
$$\Delta\left(h;f_k^0\right)\geq\Delta\left(h_k^0;f_k^0\right)\geq \Delta\left(h_k^0;f\right)   \hspace{1cm}   \forall f \in D,\; \forall h \in H_D^k$$
hold true if $h_k^0=h_k(f_k^0)$ and $h_k(f_k^0)\in H_D^k$ where $f_k^0$ are solutions to the constrained optimization problems
\begin{equation} \label{7_1}
\tilde{\Delta}_k(f)=-\Delta\left(h_k^0;f\right)=-\frac{1}{2 \pi} \int_{- \pi}^{ \pi}\left|C_{k}^{0}(e^{i \lambda})\right|^{2}\frac{f( \lambda )}{\left(f_k^{0}( \lambda )\right)^{2}}d \lambda\rightarrow \inf, \; f(\lambda) \in D,
\end{equation}
where
\[C_{1}^{0}(e^{i \lambda})=C_{N}^{0}(e^{i \lambda} )+C_{-M_1}^{0}(e^{i \lambda})= \sum\limits_{j=0}^{N}\left(B_{S_1}^{0-1}\vec{\bold{a}}_{1}\right)(j)e^{ij\lambda}+\sum\limits_{j=-\infty}^{-M_1-1}\left(B_{S_1}^{0-1}\vec{\bold{a}}_{1}\right)(j)e^{ij\lambda},\]
\[C_{2}^{0}(e^{i \lambda})=C_{N}^{0}(e^{i \lambda})+C_{N+M_2+}^{0}(e^{i \lambda})= \sum\limits_{j=0}^{N}\left(B_{S_2}^{0-1}\vec{\bold{a}}_{2}\right)(j)e^{ij\lambda}+\sum\limits_{j=N+M_2+1}^{\infty}\left(B_{S_2}^{0-1}\vec{\bold{a}}_{2}\right)(j)e^{ij\lambda},\]
\[C_{3}^{0}(e^{i \lambda})=C_{-M_1}^{0}(e^{i \lambda})+C_{N}^{0}(e^{i \lambda})+C_{N+M_2+}^{0}(e^{i \lambda})=\]
\[=\sum\limits_{j=-\infty}^{M_1-1}c^0(j)e^{ij\lambda}+\sum\limits_{j=0}^{N}c^0(j)e^{ij\lambda}+
\sum\limits_{j=N+M_2+1}^{\infty}c^0(j)e^{ij\lambda},\]
\[C_{4}^{0}(e^{i \lambda})=C_{N}^{0}(e^{i \lambda} )+C_{-M_1}^{0N_1}(e^{i \lambda})= \sum\limits_{j=0}^{N}\left(B_{S_4}^{0-1}\vec{\bold{a}}_{4}\right)(j)e^{ij\lambda}+\sum\limits_{j=-M_1-N_1}^{-M_1-1}\left(B_{S_4}^{0-1}\vec{\bold{a}}_{4}\right)(j)e^{ij\lambda},\]
\[C_{5}^{0}(e^{i \lambda})=C_{N}^{0}(e^{i \lambda})+C_{N+M_2+}^{0N_2}(e^{i \lambda})= \sum\limits_{j=0}^{N}\left(B_{S_5}^{0-1}\vec{\bold{a}}_{5}\right)(j)e^{ij\lambda}+\sum\limits_{j=N+M_2+1}^{N+M_2+N_2}\left(B_{S_5}^{0-1}\vec{\bold{a}}_{5}\right)(j)e^{ij\lambda},\]
\[C_{6}^{0}(e^{i \lambda})=C_{N}^{0}(e^{i \lambda})+C_{-M_1}^{0N_1}(e^{i \lambda})+C_{N+M_2+}^{0N_2}(e^{i \lambda})=\]
\[=\sum\limits_{j=0}^{N}\left(B_{S_6}^{0-1}\vec{\bold{a}}_{6}\right)(j)e^{ij\lambda}+\sum\limits_{j=-M_1-N-1}^{-M_1-1}\left(B_{S_6}^{0-1}\vec{\bold{a}}_{6}\right)(j)e^{ij\lambda}+\sum\limits_{j=N+M_2+1}^{N+M_2+N_2}\left(B_{S_6}^{0-1}\vec{\bold{a}}_{6}\right)(j)e^{ij\lambda}.\]

The constrained optimization problems (33) are equivalent to the unconstrained optimization problems
\begin{equation*}
\Delta_D^k(f)=\tilde{\Delta}_k(f)+\delta(f\left|D\right.)\rightarrow \inf,
\end{equation*}
where $\delta(f\left|D\right.)$ is the indicator function of the set  $D$. Solutions $f_k^0$ to this problems are characterized by the conditions $0 \in \partial\Delta_D^k(f_k^0)$ where $\partial\Delta_D^k(f_k^0)$ are the subdifferentials of the convex functionals $\Delta_D^k(f)$ at points $f_k^0$ respectively.
This condition makes it possible to find the least favourable spectral densities in some special classes of spectral densities $D$ [9], [22], [23].

Note, that the form of the functionals $\Delta\left(h_k^0;f\right)$ is convenient for application the Lagrange method of indefinite multipliers for finding solution to the problems (33). Making use the method of Lagrange multipliers and the form of subdifferentials of the indicator functions we describe relations that determine least favourable spectral densities in some special classes of spectral densities.

\section{Least favourable spectral densities in the class $D_0^{-}$}

Consider the problem of the optimal estimation of the functionals $A_{S_k}\xi$  which depend on the unknown values of a stationary stochastic sequence $\xi(j)$ in the case where the spectral density is from the class
\begin{equation*}
D_0^{-} = \left\{f(\lambda)\left|\frac{1}{2\pi}\int\limits_{-\pi}^{\pi} f^{-1}(\lambda)d\lambda \ge p\right.\right\},
\end{equation*}
where $p$ is a given number and sequences $a(k), k\in Z\backslash S_k$ which determine functionals $A_{S_k}\xi$ are strictly positive. To find solutions to the constrained optimization problems (33) we use the  Lagrange multipliers method. With the help of this method we get the equations
\begin{equation*}
\frac{\left|C_{N}^{0} (e^{i \lambda})+C_{-M_1}^{0}(e^{i \lambda})\right|^{2}}{(f_{1}^0(\lambda ))^{2}}=\frac{\alpha_1^{2}}{(f_{1}^0(\lambda ))^{2}},
\end{equation*}
\begin{equation*}
\frac{\left|C_{N}^{0} (e^{i \lambda})+C_{N+M_2+}^{0}(e^{i \lambda})\right|^{2}}{(f_{2}^0(\lambda ))^{2}}=\frac{\alpha_2^{2}}{(f_{2}^0(\lambda ))^{2}},
\end{equation*}
\begin{equation*}
\frac{\left|C_{-M_1}^{0}(e^{i \lambda})+C_{N}^{0} (e^{i \lambda})+C_{N+M_2+}^{0}(e^{i \lambda})\right|^{2}}{(f_{3}^0(\lambda ))^{2}}=\frac{\alpha_3^{2}}{(f_{3}^0(\lambda ))^{2}},
\end{equation*}
where $\alpha_k^{2}$ are unknown Lagrange multipliers. From this relations we find that the Fourier coefficients of the functions $(f_{k}^0( \lambda ))^{-1}$ satisfy the following equations respectively
\begin{equation} \label{99}
\left|\sum\limits_{j=0}^{N}c(j)e^{ij\lambda}+\sum\limits_{j=-\infty}^{-M_1-1}c(j)e^{ij\lambda}\right|^2=\alpha_1^{2},
\end{equation}
\begin{equation} \label{99_1}
\left|\sum\limits_{j=0}^{N}c(j)e^{ij\lambda}+\sum\limits_{j=N+M_2+1}^{\infty}c(j)e^{ij\lambda}\right|^2=\alpha_2^{2},
\end{equation}
\begin{equation} \label{99_2}
\left|\sum\limits_{j=-\infty}^{M_1-1}c(j)e^{ij\lambda}+\sum\limits_{j=0}^{N}c(j)e^{ij\lambda}+\sum\limits_{j=N+M_2+1}^{\infty}c(j)e^{ij\lambda}\right|^2=\alpha_3^{2},
\end{equation}
where $c(j), j \in Z\backslash S_1$, and $c(j), j \in Z\backslash S_2$, are components of the vectors $\vec{\bold{c}}_{1}$ and $\vec{\bold{c}}_{2}$ respectively, that satisfy the following equations
\[B_{S_1}^0\vec{\bold{c}}_{1}=\vec{\bold{a}}_{1},\; B_{S_2}^0\vec{\bold{c}}_{2}=\vec{\bold{a}}_{2},\]
where matrices $B_{S_1}^0$, $B_{S_2}^0$ are determined by the Fourier coefficients of the functions
\[(f_1^0(\lambda))^{-1}=\sum\limits_{k=-\infty}^{\infty}b_1^0(k)e^{ik\lambda},\; (f_2^0(\lambda))^{-1}=\sum\limits_{k=-\infty}^{\infty}b_2^0(k)e^{ik\lambda},\]
and $c(j), j \in Z\backslash S_3$ are components of the vectors $\vec{c}_N$, $\vec{c}_{-M_1-1}$, $\vec{c}_{N+M_2+1}$ that satisfy the following equations
\[\vec{a}_N=B_1^0\vec{c}_N+B_2^0\vec{c}_{-M_1-1}+B_5^0\vec{c}_{N+M_2+1},\;
\vec{a}_{-M_1-1}=B_3^0\vec{c}_N+B_4^0\vec{c}_{-M_1-1}+B_8^0\vec{c}_{N+M_2+1},\]
\[\vec{a}_{N+M_2+1}=B_6^0\vec{c}_N+B_9^0\vec{c}_{-M_1-1}+B_7^0\vec{c}_{N+M_2+1},\]
where matrices $B_k^0$ are determined by the Fourier coefficients of the function
\[(f_3^0(\lambda))^{-1}=\sum\limits_{k=-\infty}^{\infty}b_3^0(k)e^{ik\lambda}.\]

The Fourier coefficients $b_k^0(n)=b_k^0(-n), n\in Z\backslash S_k, k=1, 2,$ satisfy equations  (34), (35) respectively, and equations
\[B_{S_1}^0\vec{\bold{c}}_{1}=\vec{\bold{a}}_{1},\; B_{S_2}^0\vec{\bold{c}}_{2}=\vec{\bold{a}}_{2}.\]
These coefficients can be found from the equations
\[B_{S_1}^0\vec{\bold{\alpha}}_1=\vec{\bold{a}}_{1},\; B_{S_2}^0\vec{\bold{\alpha}}_2=\vec{\bold{a}}_{2},\]
where
\[\vec{\bold{\alpha}}_1=(0, 0, \ldots, \alpha_1, 0, 0, \ldots),\; \vec{\bold{\alpha}}_2=(\alpha_2, 0, 0, \ldots).\]

The last relations can be presented in the form of the system of equations
\[\alpha_1 b_1^0(n-N)=a(n), \, n\in Z\backslash S_1,\]
\[\alpha_2 b_2^0(n)=a(n), \, n\in Z\backslash S_2.\]

From the equations of the systems for $n=N$ and $n=0$ respectively we find the unknown values
\[\alpha_1=a(N)(b_1^0(0))^{-1},\; \alpha_2=a(0)(b_2^0(0))^{-1}.\]

It follows from the extremum conditions (31) and the restriction on the spectral densities from the class $D_0^{-}$  that the Fourier coefficients
$$b_k^0(0)=\frac{1}{2\pi}\int\limits_{-\pi}^{\pi}(f_k^0(\lambda))^{-1}d\lambda=p.$$
Thus,
  \begin{displaymath}
  b_1^0(n-N)=b_1^0(N-n)=\left\{ \begin{array}{ll}
  p\,a(n)(a(N))^{-1} & \textrm{if $n\in Z\backslash S_1$;} \\
  0 & \textrm{if $n\in S_1$,}
  \end{array} \right.
  \end{displaymath}
   \begin{displaymath}
  b_2^0(n)=b_2^0(-n)=\left\{ \begin{array}{ll}
  p\,a(n)(a(0))^{-1} & \textrm{if $n\in Z\backslash S_2$;} \\
  0 & \textrm{if $n\in S_2$.}
  \end{array} \right.
  \end{displaymath}
Hence, we can write the functions $(f_k^0(\lambda))^{-1}$, $k=1, 2$ in the form
\begin{equation} \label{dens_1}
(f_k^0(\lambda))^{-1}=\sum_{n=-\infty}^{-N-M_k-1}b_k^0(n) e^{in\lambda}+\sum_{n=-N}^{N}b_k^0(n) e^{in\lambda}+\sum_{n=N+M_k+1}^{\infty}b_k^0(n) e^{in\lambda}.
\end{equation}
In this case the functions $(f_k^0(\lambda))^{-1}$, $k=1, 2$ can be represented in the form [14]
\begin{equation}
(f_k^0(\lambda))^{-1}=\left|\sum\limits_{n=0}^{\infty}\gamma_{kn} e^{-in\lambda}\right|^2,\; \lambda \in \left[-\pi, \pi\right],
\end{equation}
where $\gamma_{kn}=0,\; n\in \left\{N+1,\hdots, N+M_k\right\}$.

The minimax spectral characteristics $h_1(f_1^0)$, $h_2(f_2^0)$ of the optimal linear estimates of the functionals $A_{S_1}\xi$, $A_{S_2}\xi$ can be calculated by the formulas (39), (40) respectively, where
 \[C_N(e^{i\lambda})+C_{-M_1}(e^{i\lambda})=\sum\limits_{j=0}^{N}c(j)e^{ij\lambda}+\sum\limits_{j=-\infty}^{-N-M_1-1}c(j)e^{ij\lambda}=\alpha_N=a(N)p^{-1},\]
\[C_N(e^{i\lambda})+C_{N+M_2+}(e^{i\lambda})=\sum\limits_{j=0}^{N}c(j)e^{ij\lambda}+\sum\limits_{j=N+M_2+1}^{\infty}c(j)e^{ij\lambda}=\alpha_1=a(0)p^{-1},\]
  namely
\begin{equation} \label{spectr_1} \begin{gathered}
h_1(f_1^0)=\sum\limits_{j=0}^{N}a(j)e^{ij\lambda}+\sum\limits_{j=-\infty}^{-M_1-1}a(j)e^{ij\lambda} - \left(\sum\limits_{k=-\infty}^{\infty}b_1^0(k) e^{ik\lambda}\right)e^{iN\lambda}a(N)p^{-1}=\\
 =-\sum\limits_{j=1}^{N}a(N-j)e^{i(N+j)\lambda}-\sum\limits_{j=N+M_1+1}^{\infty}a(N-j)e^{i(N+j)\lambda},
\end{gathered}\end{equation}
\begin{equation} \label{spectr_2} \begin{gathered}
h_2(f_2^0)=\sum\limits_{j=0}^{N}a(j)e^{ij\lambda}+\sum\limits_{j=N+M_2+1}^{\infty}a(j)e^{ij\lambda} - \left(\sum\limits_{k=-\infty}^{\infty}b_2^0(k) e^{ik\lambda}\right)a(0)p^{-1}=\\
 =-\sum\limits_{j=1}^{N}a(j)e^{-ij\lambda}-\sum\limits_{j=N+M_2+1}^{\infty}a(j)e^{-ij\lambda}.
\end{gathered}\end{equation}

Summing up our reasoning we come to conclusion that the following theorem holds true.
\begin{thm}
The least favourable in the class $D_0^{-}$ spectral densities for the optimal linear estimation of the functionals $A_{S_1}\xi$ and $A_{S_2}\xi$ determined by strictly positive sequences $a(k), k \in Z\backslash S_1$, and $a(k), k \in Z\backslash S_2$ are the spectral densities (37) with the Fourier coefficients
 \[b_1^0(n-N)=b_1^0(N-n)=p\,a(n)(a(N))^{-1}, \; n\in Z\backslash S_1,\]
 \[b_2^0(n)=b_2^0(-n)=p\,a(n)(a(0))^{-1}, \; n\in Z\backslash S_2.\]
The minimax spectral characteristics $h_1(f_1^0)$, $h_2(f_2^0)$ are given by formulas (39), (40).
The least favourable in the class $D_0^{-}$ spectral density for the optimal linear estimation of the functionals $A_{S_3}\xi$ determined by strictly positive sequences $a(k), k \in Z\backslash S_3$ satisfies relation (36)  and the optimization problem (32). The minimax spectral characteristics $h_3(f_3^0)$ can be calculated by formula (8).
\end{thm}

Consider the problem of the optimal estimation of the functionals $A_{S_k}\xi$ for $k=4, 5, 6$. To find solutions to the constrained optimization problems (33) we get the equations
\begin{equation*}
\frac{\left|C_{N}^{0} (e^{i \lambda})+C_{-M_1}^{0N_1}(e^{i \lambda})\right|^{2}}{(f_{4}^0(\lambda ))^{2}}=\frac{\alpha_4^{2}}{(f_{4}^0(\lambda ))^{2}},
\end{equation*}
\begin{equation*}
\frac{\left|C_{N}^{0} (e^{i \lambda})+C_{N+M_2+}^{0N_2}(e^{i \lambda})\right|^{2}}{(f_{5}^0(\lambda ))^{2}}=\frac{\alpha_5^{2}}{(f_{5}^0(\lambda ))^{2}},
\end{equation*}
\begin{equation*}
\frac{\left|C_{-M_1}^{0N_1}(e^{i \lambda})+C_{N}^{0} (e^{i \lambda})+C_{N+M_2+}^{0N_2}(e^{i \lambda})\right|^{2}}{(f_{6}^0(\lambda ))^{2}}=\frac{\alpha_6^{2}}{(f_{6}^0(\lambda ))^{2}},
\end{equation*}
where $\alpha_k^{2}$ are unknown Lagrange multipliers. From this relations we find that the Fourier coefficients of the functions $(f_{k}^0( \lambda ))^{-1}$ satisfy the following equations respectively
\begin{equation} \label{99_4}
\left|\sum\limits_{j=0}^{N}c(j)e^{ij\lambda}+\sum\limits_{j=-M_1-N_1}^{-M_1-1}c(j)e^{ij\lambda}\right|^2=\alpha_4^{2},
\end{equation}
\begin{equation} \label{99_5}
\left|\sum\limits_{j=0}^{N}c(j)e^{ij\lambda}+\sum\limits_{j=N+M_2+1}^{N+M_2+N_2}c(j)e^{ij\lambda}\right|^2=\alpha_5^{2},
\end{equation}
\begin{equation} \label{99_6}
\left|\sum\limits_{j=-M_1-N_1}^{M_1-1}c(j)e^{ij\lambda}+\sum\limits_{j=0}^{N}c(j)e^{ij\lambda}+\sum\limits_{j=N+M_2+1}^{N+M_2+N_2}c(j)e^{ij\lambda}\right|^2=\alpha_6^{2},
\end{equation}
where $c(j), j \in Z\backslash S_k$ are components of the vectors $\vec{\bold{c}}_{k}$
\[\vec{\bold{c}}_{4}^\top=(\vec{c}_N^\top,\vec{c}_{-M_1-1}^{N_1\top}),\; \vec{\bold{c}}_{5}^\top=(\vec{c}_N^\top,\vec{c}_{N+M_2+1}^{N_2\top}),\; \vec{\bold{c}}_{6}^\top=(\vec{c}_N^\top, \vec{c}_{-M_1-1}^{N_1\top},\vec{c}_{N+M_2+1}^{N_2\top}),\]
\[\vec{c}_N^\top=(c(0),c(1),\ldots,c(N)),\; \vec{c}_{-M_1-1}^{N_1\top}=(c(-M_1-1),c(-M_1-2),\ldots,c(-M_1-N_1)),\]
\[\vec{c}_{N+M_2+1}^{N_2\top}=(c(N+M_2+1),c(N+M_2+2),\ldots,c(N+M_2+N_2)).\]
that satisfy the equations $B_{S_k}^0\vec{\bold{c}}_{k}=\vec{\bold{a}}_{k}$ respectively where matrices $B_{S_k}^0$ are determined by the Fourier coefficients of the functions $(f_{k}^0( \lambda ))^{-1}$.

The Fourier coefficients
\[b_k^0(n)=b_k^0(-n), n\in Z\backslash S_k, k=4, 5, 6\]
satisfy equations  (41), (42), (43) respectively, and equations $B_{S_k}^0\vec{\bold{c}}_{k}=\vec{\bold{a}}_{k}$. These coefficients can be found from the equations $B_{S_k}^0\vec{\bold{\alpha}}_k=\vec{\bold{a}}_{k}$ where
\[\vec{\bold{\alpha}}_4=(0, 0, \ldots, \alpha_4, 0, 0, \ldots, 0),\; \vec{\bold{\alpha}}_5=(\alpha_5, 0, 0, \ldots, 0),\; \vec{\bold{\alpha}}_6=(0, 0, \ldots, \alpha_6).\]

The last relations can be presented in the form of the system of equations
\[\alpha_4 b_4^0(n-N)=a(n), \, n\in Z\backslash S_4,\]
\[\alpha_5 b_5^0(n)=a(n), \, n\in Z\backslash S_5,\]
\[\alpha_6 b_6^0(n-N-M_2-N_2)=a(n), \, n\in Z\backslash S_6.\]

From the equations of the systems for $n=N$, $n=0$ and $n=N+M_2+N_2$ respectively we find the unknown values \[\alpha_4=a(N)(b_4^0(0))^{-1},\;\alpha_5=a(0)(b_5^0(0))^{-1}, \alpha_6=a(N+M_2+N_2)(b_6^0(0))^{-1}.\]

It follows from the extremum conditions (31) and the restriction on the spectral densities from the class $D_0^{-}$  that the Fourier coefficients
$$b_k^0(0)=\frac{1}{2\pi}\int\limits_{-\pi}^{\pi}(f_k^0(\lambda))^{-1}d\lambda=p.$$
Thus,
  \begin{displaymath}
  b_4^0(n-N)=b_4^0(N-n)=\left\{ \begin{array}{ll}
  p\,a(n)(a(N))^{-1} & \textrm{if $n\in Z\backslash S_4$;} \\
  0 & \textrm{if $n\in S_4$,}
  \end{array} \right.
  \end{displaymath}
   \begin{displaymath}
  b_5^0(n)=b_5^0(-n)=\left\{ \begin{array}{ll}
  p\,a(n)(a(0))^{-1} & \textrm{if $n\in Z\backslash S_5$;} \\
  0 & \textrm{if $n\in S_5$,}
  \end{array} \right.
  \end{displaymath}
  \begin{displaymath}
  b_6^0(n-N-M_2-N_2)=b_6^0(N+M_2+N_2-n)=\left\{ \begin{array}{ll}
  p\,a(n)(a(N+M_2+N_2))^{-1} & \textrm{if $n\in Z\backslash S_6$;} \\
  0 & \textrm{if $n\in S_6$.}
  \end{array} \right.
  \end{displaymath}
Hence, we can write the functions $(f_k^0(\lambda))^{-1}$ in the form
\begin{equation} \label{dens_4}
(f_4^0(\lambda))^{-1}=\sum_{k=-N-M_1-N_1}^{-N-M_1-1}b_4^0(k) e^{ik\lambda}+\sum_{k=-N}^{N}b_4^0(k) e^{ik\lambda}+\sum_{k=N+M_1+1}^{N+M_1+N_1}b_4^0(k) e^{ik\lambda},
\end{equation}
\begin{equation} \label{dens_5}
(f_5^0(\lambda))^{-1}=\sum_{k=-N-M_2-N_2}^{-N-M_2-1}b_5^0(k) e^{ik\lambda}+\sum_{k=-N}^{N}b_5^0(k) e^{ik\lambda}+\sum_{k=N+M_2+1}^{N+M_2+N_2}b_5^0(k) e^{ik\lambda},
\end{equation}
\[(f_6^0(\lambda))^{-1}=\sum_{k=-N-M_2-N_2-M_1-N_1}^{-N-M_2-N_2-M_1-1}b_6^0(k) e^{ik\lambda}+\sum_{k=-N-M_2-N_2}^{-M_2-N_2}b_6^0(k) e^{ik\lambda}+\sum_{k=-(N_2-1)}^{N_2-1}b_6^0(k) e^{ik\lambda}+\]
\begin{equation} \label{dens_6}
+\sum_{k=M_2+N_2}^{N+M_2+N_2}b_6^0(k) e^{ik\lambda}+\sum_{k=N+M_2+N_2+M_1+1}^{N+M_2+N_2+M_1+N_1}b_6^0(k) e^{ik\lambda}.
\end{equation}

In this case the functions $(f_k^0(\lambda))^{-1}$, $k=4, 5, 6$ can be represented in the form [14]
\begin{equation}
(f_k^0(\lambda))^{-1}=\left|\sum\limits_{n=0}^{\infty}\gamma_{kn} e^{-in\lambda}\right|^2,\; \lambda \in \left[-\pi, \pi\right],
\end{equation}
where
\[\gamma_{4n}=0,\; N+1\le n\le N+M_1,\; n> N+M_1+N_1,\]
\[\gamma_{5n}=0,\; N+1\le n\le N+M_2,\; n> N+M_2+N_2,\]
\[\gamma_{6n}=0, N_2\le n< M_2+N_2, N+M_2+N_2<n\le N+M_2+N_2+M_1, n> N+M_2+N_2+M_1+N_1.\]

The minimax spectral characteristics $h_k(f_k^0)$ of the optimal linear estimates of the functionals $A_{S_k}\xi$ can be calculated by the formulas
\begin{equation} \label{spectr_4}
h_4(f_4^0)=-\sum\limits_{j=1}^{N}a(N-j)e^{i(N+j)\lambda}-\sum\limits_{j=N+M_1+1}^{N+M_1+N_1}a(N-j)e^{i(N+j)\lambda},
\end{equation}
\begin{equation} \label{spectr_5}
h_5(f_5^0)=-\sum\limits_{j=1}^{N}a(j)e^{-ij\lambda}-\sum\limits_{j=N+M_2+1}^{N+M_2+N_2}a(j)e^{-ij\lambda},
\end{equation}
\[h_6(f_6^0)=-\sum\limits_{j=1}^{N_2-1}a((N+M_2+N_2-j)e^{i((N+M_2+N_2+j)\lambda}-\sum\limits_{j=M_2+N_2}^{N+M_2+N_2}a(N+M_2+N_2-j)e^{i(N+M_2+N_2+j)\lambda}-\]
\begin{equation} \label{spectr_6}
-\sum\limits_{j=N+M_2+N_2+M_1+1}^{N+M_2+N_2+M_1+N_1}a(N+M_2+N_2-j)e^{i(N+M_2+N_2+j)\lambda}.
\end{equation}

Summing up our reasoning we come to conclusion that the following corollary holds true.
\begin{nas}
The least favourable in the class $D_0^{-}$ spectral densities for the optimal linear estimation of the functionals $A_{S_k}\xi$, $k=4, 5, 6$ determined by strictly positive sequences $a(k), k \in Z\backslash S_k$ are the spectral densities (44), (45), (46) with the Fourier coefficients
\[b_4^0(n-N)=b_4^0(N-n)=p\,a(n)(a(N))^{-1}, \; n\in Z\backslash S_4,\]
\[b_5^0(n)=b_5^0(-n)=p\,a(n)(a(0))^{-1}, \; n\in Z\backslash S_5,\]
\[ b_6^0(n-N-M_2-N_2)=b_6^0(N+M_2+N_2-n)=p\,a(n)(a(N+M_2+N_2))^{-1}, \; n\in Z\backslash S_6.\]
The minimax spectral characteristics $h_k(f_k^0)$ are given by formulas (48), (49), (50).
\end{nas}

\section{Least favourable spectral densities in the class $D_W$}

Consider the problem of the optimal estimation of the functionals $A_{S_k}\xi$ which depend on the unknown values of a stationary stochastic sequence $\xi(j)$ in the case where the spectral density of  the  sequence is from the set of spectral densities
\begin{equation*}
D_W = \left\{f(\lambda)\left|\frac{1}{2\pi}\int\limits_{-\pi}^{\pi} (f(\lambda))^{-1}\cos(n\lambda) d\lambda= b(n),\; n=0, 1, \ldots, W  \right. \right\},
\end{equation*}
where $b(n),\; n=0, 1, \ldots, W$ is a strictly positive sequence. There is an infinite number of functions in the class $D_W$ [14] and the function
\[(f(\lambda))^{-1}=\sum_{k=-W}^{W}b(|n|) e^{ik\lambda}>0, \; \lambda \in \left[-\pi,\pi\right].\]
To find solutions to the constrained optimization problems (33) we use the  Lagrange multipliers method. With the help of this method we get the equations
\begin{equation}\label{w_1}
\left|\sum\limits_{j=0}^{N}c(j)e^{ij\lambda}+\sum\limits_{j=-\infty}^{-M_1-1}c(j)e^{ij\lambda}\right|^2=\sum_{n=0}^{W}\alpha_1(n)\cos(n\lambda)=\left|\sum_{n=0}^{W}p_1(n)e^{in\lambda}\right|,
\end{equation}
\begin{equation}\label{w_2}
\left|\sum\limits_{j=0}^{N}c(j)e^{ij\lambda}+\sum\limits_{j=N+M_2+1}^{\infty}c(j)e^{ij\lambda}\right|^2=\sum_{n=0}^{W}\alpha_2(n)\cos(n\lambda)=\left|\sum_{n=0}^{W}p_2(n)e^{in\lambda}\right|,
\end{equation}
\begin{equation}\label{w_3}
\left|\sum\limits_{j=-\infty}^{-M_1-1}c(j)e^{ij\lambda}+\sum\limits_{j=0}^{N}c(j)e^{ij\lambda}+\sum\limits_{j=N+M_2+1}^{\infty}c(j)e^{ij\lambda}\right|^2=\sum_{n=0}^{W}\alpha_3(n)\cos(n\lambda)=\left|\sum_{n=0}^{W}p_3(n)e^{in\lambda}\right|,
\end{equation}
where $c(j), j \in Z\backslash S_1$, and $c(j), j \in Z\backslash S_2$, are components of the vectors $\vec{\bold{c}}_{1}$ and $\vec{\bold{c}}_{2}$ respectively, that satisfy the following equations
\[B_{S_1}^0\vec{\bold{c}}_{1}=\vec{\bold{a}}_{1},\; B_{S_2}^0\vec{\bold{c}}_{2}=\vec{\bold{a}}_{2},\]
and $c(j), j \in Z\backslash S_3$ are components of the vectors $\vec{c}_N$, $\vec{c}_{-M_1-1}$, $\vec{c}_{N+M_2+1}$ that satisfy the following equations
\[\vec{a}_N=B_1^0\vec{c}_N+B_2^0\vec{c}_{-M_1-1}+B_5^0\vec{c}_{N+M_2+1},\;
\vec{a}_{-M_1-1}=B_3^0\vec{c}_N+B_4^0\vec{c}_{-M_1-1}+B_8^0\vec{c}_{N+M_2+1},\]
\[\vec{a}_{N+M_2+1}=B_6^0\vec{c}_N+B_9^0\vec{c}_{-M_1-1}+B_7^0\vec{c}_{N+M_2+1}.\]

Let $M_1\ge N$. Then matrices $B_{S_1}^0$, $B_{S_2}^0$ are determined by the known
\[b_k^0(n)=b(|n|),\; n\in Z\backslash S_k \; \text{and}\; |n|\in\{0, 1, \ldots, W\}\]
and the unknown Fourier coefficients
\[b_k^0(n),\; n\in Z\backslash S_k\; \text{and}\; |n|\notin\{0, 1, \ldots, W\}\]
of the functions $(f_1^0(\lambda))^{-1}$, $(f_2^0(\lambda))^{-1}$. The unknown coefficients
\[p_k(n),\; n\in Z\backslash S_k\; \text{and}\; |n|\notin\{0, 1, \ldots, W\},\]
and $b_k^0(n)$ can be found from equations
\[B_{S_1}^0\vec{p}_{1}^0=\vec{\bold{a}}_{1},\; B_{S_2}^0\vec{p}_{2}^0=\vec{\bold{a}}_{2},\]
with \[\vec{p}_{k}^0=(p_k(0), \ldots, p_k(W_k),0, 0, \ldots),\]
where
\[W_1=W\; \text{if}\; W\le N,\; W_1=N\; \text{if}\; N<W\le M_1,\; W_1=N+W-M_1\; \text{if}\; W>M_1,\]
\[W_2=W\; \text{if}\; W\le N,\; W_2=N\; \text{if}\; N<W\le N+M_2,\; W_2=W-M_2\; \text{if}\; W>N+M_2.\]

If the sequences $b_k^0(n)$ that are constructed from the strictly positive sequence $b(n),\; n=0, 1, \ldots, W$ and calculated coefficients $b_k^0(n)$, $n\in Z\backslash S_k$ and $|n|\notin\{0, 1, \ldots, W\}$ are also strictly positive then the least favourable spectral densities $f_k^0(\lambda)$, $k=1, 2$ are determined by Fourier coefficients $b_k^0(n)$ of the functions
$(f_k^0(\lambda))^{-1}$, $k=1, 2$
\begin{equation}\label{w_4}
(f_k^0(\lambda))^{-1}=\sum_{n=-W}^{W}(b(|n|) e^{in\lambda}+\sum_{n\in Z\backslash S_k, \; |n|\notin\{0, 1, \ldots, W\}}(b_k^0(n) e^{in\lambda}+b_k^0(n) e^{-in\lambda})=\left|\sum\limits_{n=0}^{\infty}\gamma_{kn} e^{-in\lambda}\right|^2,
\end{equation}
where $\gamma_{1n}=0,\; n\in \left\{N+1,\hdots, M_1\right\}\backslash\{0, 1, \ldots, W\}$, $\gamma_{2n}=0,\; n\in \left\{N+1,\hdots, N+M_2\right\}\backslash\{0, 1, \ldots, W\}$.

Summing up our reasoning we come to conclusion that the following theorem holds true.
\begin{thm}
The least favourable in the class $D_W$ spectral densities $f_k^0(\lambda)$, $k=1, 2$ for the optimal linear estimation of the functionals $A_{S_1}\xi$ and $A_{S_2}\xi$ in the case where the strictly positive sequence $b(n),\; n=0, 1, \ldots, W$ and solutions $b_k^0(n)$, $n\in Z\backslash S_k$ and $|n|\notin\{0, 1, \ldots, W\}$ of the equations $B_{S_k}^0\vec{p}_{k}^0=\vec{\bold{a}}_{k}$ form a strictly positive sequence are determined by equations (54) . The minimax spectral characteristics $h_1(f_1^0)$, $h_2(f_2^0)$ are given by formulas (6), (7). The least favourable in the class $D_W$ spectral density for the optimal linear estimation of the functionals $A_{S_3}\xi$ satisfies relation (53) and the optimization problem (32). The minimax spectral characteristics $h_3(f_3^0)$ can be calculated by formula (32).
\end{thm}

Consider the problem of the optimal estimation of the functionals $A_{S_k}\xi$ for $k=4, 5, 6$. To find solutions to the constrained optimization problems (33) we get the equations
\begin{equation}\label{w_5}
\left|\sum\limits_{j=0}^{N}c(j)e^{ij\lambda}+\sum\limits_{j=-M_1-N_1}^{-M_1-1}c(j)e^{ij\lambda}\right|^2=\sum_{n=0}^{W}\alpha_4(n)\cos(n\lambda)=\left|\sum_{n=0}^{W}p_4(n)e^{in\lambda}\right|,
\end{equation}
\begin{equation}\label{w_6}
\left|\sum\limits_{j=0}^{N}c(j)e^{ij\lambda}+\sum\limits_{j=N+M_2+1}^{N+M_2+N_2}c(j)e^{ij\lambda}\right|^2=\sum_{n=0}^{W}\alpha_5(n)\cos(n\lambda)=\left|\sum_{n=0}^{W}p_5(n)e^{in\lambda}\right|,
\end{equation}
\begin{equation}\label{w_7}
\left|\sum\limits_{j=-M_1-N_1}^{-M_1-1}c(j)e^{ij\lambda}+\sum\limits_{j=0}^{N}c(j)e^{ij\lambda}+\sum\limits_{j=N+M_2+1}^{N+M_2+N_2}c(j)e^{ij\lambda}\right|^2=\sum_{n=0}^{W}\alpha_6(n)\cos(n\lambda)=\left|\sum_{n=0}^{W}p_6(n)e^{in\lambda}\right|,
\end{equation}
where $c(j), j \in Z\backslash S_k$ are components of the vectors $\vec{\bold{c}}_{k}$ that satisfy equations $B_{S_k}^0\vec{\bold{c}}_{k}=\vec{\bold{a}}_{k}$.

Let $M_1\ge N$, $N+M_2\ge M_1+N_1$. Let, first, the functional $A_{S_4}\xi$ satisfy condition $W\ge M_1+N_1$ and the functionals $A_{S_5}\xi$, $A_{S_6}\xi$ satisfy condition $W\ge N+M_2+N_2$. In this case the given Fourier coefficients $b(n),\; n=0, 1, \ldots, W$ define matrices $B_{S_k}^0$, $k=4, 5, 6$ and the optimization problem (31) is degenerated. Thus every density $f(\lambda)\in D_W$ is the least favourable and the density of the autoregressive stochastic sequence
\begin{equation}\label{w_8}
(f^0(\lambda))^{-1}=\sum_{n=-W}^{W}b(|n|) e^{in\lambda}=\left|\sum\limits_{n=0}^{W}\gamma_{n} e^{-in\lambda}\right|^2,
\end{equation}
is the least favourable, as well.

Let the functional $A_{S_4}\xi$ satisfy condition $W<M_1+N_1$ and the functionals $A_{S_5}\xi$, $A_{S_6}\xi$ satisfy condition $W<N+M_2+N_2$.Then matrices $B_{S_k}^0$, $k=4, 5, 6$ are determined by the known
\[b_k^0(n)=b(|n|),\; n\in Z\backslash S_k \; \text{and}\; |n|\in\{0, 1, \ldots, W\}\]
and the unknown Fourier coefficients
\[b_k^0(n),\; n\in Z\backslash S_k\; \text{and}\; |n|\notin\{0, 1, \ldots, W\}\]
of the functions $(f_k^0(\lambda))^{-1}$, $k=4, 5, 6$. The unknown coefficients
\[p_k(n),\; n\in Z\backslash S_k\; \text{and}\; |n|\notin\{0, 1, \ldots, W\}\]
and $b_k^0(n)$ can be found from equations $B_{S_k}^0\vec{p}_{k}^0=\vec{\bold{a}}_{k}$ with $\vec{p}_{k}^0=(p_k(0), \ldots, p_k(W_k),0, 0, \ldots)$ where
\[W_4=W\; \text{if}\; W\le N,\; W_4=N\; \text{if}\; N<W\le M_1,\; W_4=N+W-M_1\; \text{if}\; M_1<W<M_1+N_1,\]
\[W_5=W\; \text{if}\; W\le N,\; W_5=N\; \text{if}\; N<W\le N+M_2,\; W_5=W-M_2\; \text{if}\; N+M_2<W<N+M_2+N_2,\]
\[W_6=W\; \text{if}\; W\le N,\; W_6=N\; \text{if}\; N<W\le M_1,\; W_6=N+W-M_1\; \text{if}\; M_1<W\le M_1+N_1,\]
\[W_6=N+N_1\; \text{if}\; M_1+N_1<W\le N+M_2,\; W_6=N_1+W-M_2\; \text{if}\; N+M_2<W<N+M_2+N_2.\]
If the sequences $b_k^0(n)$ that are constructed from the strictly positive sequence $b(n),\; n=0, 1, \ldots, W$ and calculated coefficients $b_k^0(n)$, $n\in Z\backslash S_k$ and $|n|\notin\{0, 1, \ldots, W\}$ are also strictly positive then the least favourable spectral densities $f_k^0(\lambda)$, $k=4, 5, 6$ are determined by Fourier coefficients $b_k^0(n)$ of the functions $(f_k^0(\lambda))^{-1}$, $k=4, 5, 6$
\begin{equation}\label{w_9}
(f_k^0(\lambda))^{-1}=\sum_{n=-W}^{W}b(|n|) e^{in\lambda}+\sum_{n\in Z\backslash S_k, \; |n|\notin\{0, 1, \ldots, W\}}(b_k^0(n) e^{in\lambda}+b_k^0(n) e^{-in\lambda})=\left|\sum\limits_{n=0}^{\infty}\gamma_{kn} e^{-in\lambda}\right|^2,
\end{equation}
where
\[\gamma_{4n}=0,\; n\in \left\{N+1,\hdots, M_1\right\}\backslash\{0, 1, \ldots, W\},\; n> M_1+N_1,\]
\[\gamma_{5n}=0,\; n\in \left\{N+1,\hdots, N+M_2\right\}\backslash\{0, 1, \ldots, W\},\; n> N+M_2+N_2,\]
\[\gamma_{6n}=0, n\in \left\{N+1,\hdots, M_1\right\}\cup\left\{M_1+N_1+1,\hdots, N+M_2\right\}\backslash\{0, 1, \ldots, W\},\; n> N+M_2+N_2.\]

Summing up our reasoning we come to conclusion that the following corollary holds true.
\begin{nas}
The least favourable in the class $D_W$ spectral densities $f_k^0(\lambda)$, $k=4, 5, 6$ for the optimal linear estimation of the functionals $A_{S_k}\xi$, $k=4, 5, 6$ 
in the case where the functional $A_{S_4}\xi$ satisfy condition $W\ge M_1+N_1$ and the functionals $A_{S_5}\xi$, 
$A_{S_6}\xi$ satisfy condition $W\ge N+M_2+N_2$ are the spectral densities (58) of autoregressive stochastic sequence of order $W$
 determined by coefficients $b(n),\; n=0, 1, \ldots, W$. In the case where the functional $A_{S_4}\xi$ satisfy condition $W<M_1+N_1$ and the functionals $A_{S_5}\xi$, $A_{S_6}\xi$ 
 satisfy condition $W<N+M_2+N_2$ and the strictly positive sequence $b(n),\; n=0, 1, \ldots, W$ and solutions $b_k^0(n)$, $n\in Z\backslash S_k$ and $|n|\notin\{0, 1, \ldots, W\}$ 
 of the equations $B_{S_k}^0\vec{p}_{k}^0=\vec{\bold{a}}_{k}$ form a strictly positive sequence, 
 the least favourable in the class $D_W$ spectral densities $f_k^0(\lambda)$, $k=4, 5, 6$ are determined by equations (59). The minimax spectral characteristics $h_k(f_k^0)$ are given by formulas (6), (7), (8).
\end{nas}

\section{Least favourable spectral densities in the class $D_v^{u}$}

Consider the problem of the optimal estimation of the functionals $A_{S_k}\xi$ which depend on the unknown values of a stationary stochastic sequence $\xi(j)$ in the case where the spectral density of  the  sequence is from the set of spectral densities
\begin{equation*}
D_v^u = \left\{f(\lambda)\left|v(\lambda)\leq f(\lambda)\leq u(\lambda),\; \frac{1}{2\pi}\int\limits_{-\pi}^{\pi} (f(\lambda))^{-1} d\lambda= p \right. \right\},
\end{equation*}
where $v(\lambda), u(\lambda)$ are given spectral densities. To find solutions to the constrained optimization problems (33) for the set $D_v^u$ of admissible spectral densities we use the conditions $0 \in \partial\Delta_D^k(f_k^0)$.
It follows from the conditions $0 \in \partial\Delta_D^k(f_k^0)$ for $D=D_v^u$ that the Fourier coefficients of the functions $(f_k^0(\lambda))^{-1}$, $k=1, 2, 3$  satisfy the equations
\begin{equation}\label{19}
\left|\sum\limits_{j=0}^{N}\left(B_{S_1}^{0-1}\vec{\bold{a}}_1\right)(j)e^{ij\lambda}+\sum\limits_{j=-\infty}^{-M_1-1}\left(B_{S_1}^{-1}\vec{\bold{a}}_1\right)(j)e^{ij\lambda}\right|^2=\gamma_{11}(\lambda)+\gamma_{12}(\lambda)+\alpha_1^2,
\end{equation}
\begin{equation}\label{19_1}
\left|\sum\limits_{j=0}^{N}\left(B_{S_2}^{0-1}\vec{\bold{a}}_2\right)(j)e^{ij\lambda}+\sum\limits_{j=N+M_2+1}^{\infty}\left(B_{S_2}^{-1}\vec{\bold{a}}_2\right)(j)e^{ij\lambda}\right|^2=\gamma_{21}(\lambda)+\gamma_{22}(\lambda)+\alpha_2^2,
\end{equation}
\begin{equation}\label{19_2}
\left|C_{-M_1}^{0}(e^{i \lambda})+C_{N}^{0} (e^{i \lambda})+C_{N+M_2+}^{0}(e^{i \lambda})\right|^{2}=\gamma_{31}(\lambda)+\gamma_{32}(\lambda)+\alpha_3^{2},
\end{equation}
where $\gamma_{k1}(\lambda)\geq 0$ and $\gamma_{k1}(\lambda)=0$ if $f_k^0(\lambda)\geq v(\lambda);$ $\gamma_{k2}(\lambda)\leq 0$ and $\gamma_{k2}(\lambda)=0$ if $f_k^0(\lambda)\leq u(\lambda).$
Therefore, in the case where  $v(\lambda)\leq f_k^0(\lambda)\leq u(\lambda)$, the functions $(f_k^0(\lambda))^{-1}, \; k=1, 2$ are of the form
\begin{equation} \label{uv}
(f_k^0(\lambda))^{-1}=\sum_{n=-\infty}^{-N-M_k-1}b_k^0(n) e^{in\lambda}+\sum_{n=-N}^{N}b_k^0(n) e^{in\lambda}+\sum_{n=N+M_k+1}^{\infty}b_k^0(n) e^{in\lambda}.
\end{equation}
where
\[b_1^0(n-N)=b_1^0(N-n)=p\,a(n)(a(N))^{-1}, \; n\in Z\backslash S_1,\]
\[b_2^0(n)=b_2^0(-n)=p\,a(n)(a(0))^{-1}, \; n\in Z\backslash S_2.\]
The densities (63) are the least favourable in the class $D_v^u$ if the following inequalities holds true
\begin{equation}\label{18}
v(\lambda) \leq \sum_{n=-\infty}^{-N-M_k-1}b_k^0(n) e^{in\lambda}+\sum_{n=-N}^{N}b_k^0(n) e^{in\lambda}+\sum_{n=N+M_k+1}^{\infty}b_k^0(n) e^{in\lambda}\leq u(\lambda), \,\, \lambda \in \left[-\pi,\pi\right].
\end{equation}

The following theorem holds true.
\begin{thm}
If the coefficients
\[b_1^0(n-N)=b_1^0(N-n)=p\,a(n)(a(N))^{-1}, \; n\in Z\backslash S_1,\]
\[b_2^0(n)=b_2^0(-n)=p\,a(n)(a(0))^{-1}, \; n\in Z\backslash S_2,\]
satisfy the inequalities (64) then the least favourable in the class $D_v^{u}$ spectral densities $f_k^0(\lambda)$ 
for the optimal linear estimate of the functionals $A_{S_k}\xi, \; k=1, 2$ are determined by formulas (63). 
The minimax spectral characteristics $h_k(f_k^0)$ of the estimates can be calculated by the formulas (39), (40). 
If the inequalities (64) are not satisfied, then the least favourable spectral densities $f_k^0(\lambda)$, $k=1, 2, 3$ 
in the class $D_v^{u}$ for the optimal linear estimate of the functionals $A_{S_k}\xi, \; k=1, 2, 3$ 
are determined by relations (60), (61), (62) and the optimization problems (31), (32). 
The spectral minimax characteristics $h_k(f_k^0)$ of the estimates can be calculated by formulas (6), (7), (8).
\end{thm}

Consider the problem of the optimal estimation of the functionals $A_{S_k}\xi$ for $k=4, 5, 6$. To find solutions to the constrained optimization problems (33) we get the equations
\begin{equation}\label{19_4}
\left|\sum\limits_{j=0}^{N}\left(B_{S_4}^{0-1}\vec{\bold{a}}_4\right)(j)e^{ij\lambda}+\sum\limits_{j=-M_1-N_1}^{-M_1-1}\left(B_{S_4}^{-1}\vec{\bold{a}}_4\right)(j)e^{ij\lambda}\right|^2=\gamma_{41}(\lambda)+\gamma_{42}(\lambda)+\alpha_4^2,
\end{equation}
\begin{equation}\label{19_5}
\left|\sum\limits_{j=0}^{N}\left(B_{S_5}^{0-1}\vec{\bold{a}}_5\right)(j)e^{ij\lambda}+\sum\limits_{j=N+M_2+1}^{N+M_2+N_2}\left(B_{S_5}^{-1}\vec{\bold{a}}_5\right)(j)e^{ij\lambda}\right|^2=\gamma_{51}(\lambda)+\gamma_{52}(\lambda)+\alpha_5^2,
\end{equation}
\[\left|\sum\limits_{j=0}^{N}\left(B_{S_6}^{0-1}\vec{\bold{a}}_6\right)(j)e^{ij\lambda}+\sum\limits_{j=-M_1-N_1}^{-M_1-1}\left(B_{S_6}^{-1}\vec{\bold{a}}_6\right)(j)e^{ij\lambda}+\sum\limits_{j=N+M_2+1}^{N+M_2+N_2}\left(B_{S_6}^{-1}\vec{\bold{a}}_6\right)(j)e^{ij\lambda}\right|^2=\]
\begin{equation}\label{19_6}
=\gamma_{61}(\lambda)+\gamma_{62}(\lambda)+\alpha_6^2,
\end{equation}
where $\gamma_{k1}(\lambda)\geq 0$ and $\gamma_{k1}(\lambda)=0$ if $f_k^0(\lambda)\geq v(\lambda);$ $\gamma_{k2}(\lambda)\leq 0$ and $\gamma_{k2}(\lambda)=0$ if $f_k^0(\lambda)\leq u(\lambda).$
Therefore, in the case where  $v(\lambda)\leq f_k^0(\lambda)\leq u(\lambda)$ the functions $(f_k^0(\lambda))^{-1}, \; k=4, 5, 6$ are of the form
\begin{equation} \label{uv_4}
(f_4^0(\lambda))^{-1}=\sum_{k=-N-M_1-N_1}^{-N-M_1-1}b_4^0(k) e^{ik\lambda}+\sum_{k=-N}^{N}b_4^0(k) e^{ik\lambda}+\sum_{k=N+M_1+1}^{N+M_1+N_1}b_4^0(k) e^{ik\lambda},
\end{equation}
\begin{equation} \label{uv_5}
(f_5^0(\lambda))^{-1}=\sum_{k=-N-M_2-N_2}^{-N-M_2-1}b_5^0(k) e^{ik\lambda}+\sum_{k=-N}^{N}b_5^0(k) e^{ik\lambda}+\sum_{k=N+M_2+1}^{N+M_2+N_2}b_5^0(k) e^{ik\lambda},
\end{equation}
\[(f_6^0(\lambda))^{-1}=\sum_{k=-N-M_2-N_2-M_1-N_1}^{-N-M_2-N_2-M_1-1}b_6^0(k) e^{ik\lambda}+\sum_{k=-N-M_2-N_2}^{-M_2-N_2}b_6^0(k) e^{ik\lambda}+\sum_{k=-(N_2-1)}^{N_2-1}b_6^0(k) e^{ik\lambda}+\]
\begin{equation} \label{uv_6}
+\sum_{k=M_2+N_2}^{N+M_2+N_2}b_6^0(k) e^{ik\lambda}+\sum_{k=N+M_2+N_2+M_1+1}^{N+M_2+N_2+M_1+N_1}b_6^0(k) e^{ik\lambda},
\end{equation}
where
\[b_4^0(n-N)=b_4^0(N-n)=p\,a(n)(a(N))^{-1}, \; n\in Z\backslash S_4,\]
\[b_5^0(n)=b_5^0(-n)=p\,a(n)(a(0))^{-1}, \; n\in Z\backslash S_5,\]
\[ b_6^0(n-N-M_2-N_2)=b_6^0(N+M_2+N_2-n)=p\,a(n)(a(N+M_2+N_2))^{-1}, \; n\in Z\backslash S_6.\]
The densities (68), (69), (70) are the least favourable in the class $D_v^u$ if the following inequalities holds true
\begin{equation}\label{18_4}
v(\lambda) \leq \sum_{k=-N-M_1-N_1}^{-N-M_1-1}b_4^0(k) e^{ik\lambda}+\sum_{k=-N}^{N}b_4^0(k) e^{ik\lambda}+\sum_{k=N+M_1+1}^{N+M_1+N_1}b_4^0(k) e^{ik\lambda}\leq u(\lambda), \,\, \lambda \in \left[-\pi,\pi\right],
\end{equation}
\begin{equation}\label{18_5}
v(\lambda) \leq \sum_{k=-N-M_2-N_2}^{-N-M_2-1}b_5^0(k) e^{ik\lambda}+\sum_{k=-N}^{N}b_5^0(k) e^{ik\lambda}+\sum_{k=N+M_2+1}^{N+M_2+N_2}b_5^0(k) e^{ik\lambda}\leq u(\lambda), \,\, \lambda \in \left[-\pi,\pi\right],
\end{equation}
\[v(\lambda) \leq \sum_{k=-N-M_2-N_2-M_1-N_1}^{-N-M_2-N_2-M_1-1}b_6^0(k) e^{ik\lambda}+\sum_{k=-N-M_2-N_2}^{-M_2-N_2}b_6^0(k) e^{ik\lambda}+\sum_{k=-(N_2-1)}^{N_2-1}b_6^0(k) e^{ik\lambda}+\]
\begin{equation}\label{18_6}
+\sum_{k=M_2+N_2}^{N+M_2+N_2}b_6^0(k) e^{ik\lambda}+\sum_{k=N+M_2+N_2+M_1+1}^{N+M_2+N_2+M_1+N_1}b_6^0(k) e^{ik\lambda}\leq u(\lambda), \,\, \lambda \in \left[-\pi,\pi\right].
\end{equation}

The following corollary holds true.

\begin{nas}
If the coefficients
\[b_4^0(n-N)=b_4^0(N-n)=p\,a(n)(a(N))^{-1}, \; n\in Z\backslash S_4,\]
\[b_5^0(n)=b_5^0(-n)=p\,a(n)(a(0))^{-1}, \; n\in Z\backslash S_5,\]
\[ b_6^0(n-N-M_2-N_2)=b_6^0(N+M_2+N_2-n)=p\,a(n)(a(N+M_2+N_2))^{-1}, \; n\in Z\backslash S_6,\]
satisfy the inequalities  (71), (72), (73) then the least favourable in the class $D_v^{u}$ spectral densities $f_k^0(\lambda)$ for the optimal linear estimate of the functionals $A_{S_k}\xi, \; k=4, 5, 6$ are determined by formulas (68), (69), (70). The minimax spectral characteristics $h_k(f_k^0)$ of the estimates can be calculated by the formulas (48), (49), (50). If the inequalities  (71), (72), (73) are not satisfied, then the least favourable spectral densities $f_k^0(\lambda)$, $k=4, 5, 6$ in the class $D_v^{u}$ for the optimal linear estimate of the functionals $A_{S_k}\xi, \; k=4, 5, 6$ are determined by relations (65), (66), (67) and the optimization problems (31). The spectral minimax characteristics $h_k(f_k^0)$ of the estimates can be calculated by formulas (25), (27), (29).
\end{nas}

\section{Conclusions}

In this article we describe methods of solution of the problem of the mean-square optimal linear estimation of the functionals
\[A_{S_1}\xi=\sum\limits_{j=-\infty}^{-M_1-1}a(j)\xi(j)+\sum\limits_{j=0}^{N}a(j)\xi(j),\; A_{S_2}\xi=\sum\limits_{j=0}^{N}a(j)\xi(j)+\sum\limits_{j=N+M_2+1}^{\infty}a(j)\xi(j),\]
\[A_{S_3}\xi=\sum\limits_{j=-\infty}^{-M_1-1}a(j)\xi(j)+
\sum\limits_{j=0}^{N}a(j)\xi(j)+\sum\limits_{j=N+M_2+1}^{\infty}a(j)\xi(j),\]
which depend on the unknown values of a stochastic stationary sequence $\xi(j),\,j\in \mathbb{Z}$. Estimates are based on observations of the sequence $\xi(j)$ at points $j\in S_k$ respectively, where
\[S_1=\{-M_1, \ldots, -1\}\cup \{N+1, N+2, \ldots\}, \; S_2=\{\ldots, -2, -1\}\cup  \{N+1, \ldots, N+M_2\},\]
\[S_3=\{-M_1, \ldots, -1\}\cup \{N+1, \ldots, N+M_2\}.\]
We provide formulas for calculating the values of the mean square errors and the spectral characteristics of the optimal linear estimates of the functionals in the case where the spectral density of the sequence $\xi(j)$ is exactly known. In the case where the spectral density is unknown while a set of admissible spectral densities is given, the minimax approach is applied. We obtain formulas that determine the least favourable spectral densities and the minimax spectral characteristics of the optimal linear estimates of the functionals $A_{S_k}\xi$ for concrete classes of admissible spectral densities.

\end{document}